\documentclass[11pt]{amsart}

\usepackage{amsmath}
\usepackage{amssymb}
\usepackage{epsfig}

\def\S{\mathfrak S}
\def\R{\mathcal R}

\DeclareMathOperator{\aaa}{\bf{a}}
\DeclareMathOperator{\Alt}{Alt}
\DeclareMathOperator{\bbb}{\bf{b}}
\DeclareMathOperator{\ccc}{\bf{c}}
\DeclareMathOperator{\cd}{cd}
\DeclareMathOperator{\co}{co}
\DeclareMathOperator{\code}{code}

\DeclareMathOperator{\ddd}{\bf{d}}
\DeclareMathOperator{\des}{des}

\DeclareMathOperator{\Exp}{Exp}
\DeclareMathOperator{\inv}{inv}

\DeclareMathOperator{\maj}{maj}

\DeclareMathOperator{\occp}{occ_{31-2}}

\DeclareMathOperator{\SSS}{SS}

\DeclareMathOperator{\STAT}{\textsc{stat}}

\newtheorem{theorem}{Theorem}[section]
\newtheorem{corollary}[theorem]{Corollary}
\newtheorem{lemma}[theorem]{Lemma}
\newtheorem{proposition}[theorem]{Proposition}

\newtheorem{problem}[theorem]{Problem}
\newtheorem{definition}[theorem]{Definition}

\newenvironment{proof_}[1]{\noindent {\it #1}}
                         {{\qed}}

\title{Variations on Descents and Inversions in Permutations}
\author{Denis Chebikin}

\address{Department of Mathematics, Massachusetts Institute of Technology, Cambridge, MA 02139}
\email{chebikin@gmail.com}

\begin{document}

\begin{abstract}
We study new statistics on permutations that are variations on the descent and
the inversion statistics. 
In particular, we consider the \emph{alternating descent set} of a permutation
$\sigma = \sigma_1\sigma_2\cdots\sigma_n$ defined as the set of indices $i$
such that either $i$ is odd and $\sigma_i > \sigma_{i+1}$, or $i$ is even and
$\sigma_i < \sigma_{i+1}$. We show that this statistic is equidistributed with
the \emph{$3$-descent set} statistic on permutations $\tilde{\sigma}
= \sigma_1\sigma_2\cdots\sigma_{n+1}$ with $\sigma_1=1$, defined to be the set
of indices $i$ such that the triple $\sigma_i \sigma_{i+1} \sigma_{i+2}$ forms
an odd permutation of size $3$. We then introduce Mahonian inversion statistics
corresponding to the two new variations of descents and show that the joint
distributions of the resulting descent-inversion pairs are the same.
We examine the generating functions involving \emph{alternating Eulerian
polynomials}, defined by analogy with the classical Eulerian polynomials
$\sum_{\sigma\in\S_n} t^{\des(\sigma)+1}$ using alternating descents.
For the alternating descent set statistic,
we define the generating polynomial in two non-commutative variables by
analogy with the $ab$-index of the Boolean algebra $B_n$, and make
observations about it.
By
looking at the number of alternating inversions in alternating (down-up)
permutations, we obtain a new $q$-analog of the Euler number $E_n$ and show
how it emerges in a $q$-analog of an identity expressing $E_n$ as a weighted
sum of Dyck paths.
\end{abstract}

\maketitle

\section{Introduction}\label{intro}

Specifying the descent set of a permutation can be thought of as giving 
information on how the elements are ordered locally, namely, which pairs of
consecutive elements are ordered properly and which are not, the latter
constituting the descents. The original idea that became the starting point
of this research was to generalize descent sets to indicators of relative
orders of $k$-tuples of consecutive elements, the next simplest case being
$k=3$. In this case there are $6$ possible relative orders,
and thus the
analog of the descent set enumerator $\Psi_n(\aaa,\bbb)$,
also known as the $ab$-index of the Boolean algebra $B_n$,
would involve $6$ non-commuting variables.
In order to defer overcomplication, to keep the number of variables at $2$,
and to stay close to classical permutation statistics, we can divide triples
of consecutive elements into merely
``proper'' or ``improper'', defined as having
the relative order of an even or an odd permutation of size $3$, respectively.
We call the improper triples \emph{$3$-descents}, and denote
the set of positions at which $3$-descents occur in a permutation $\sigma$
by $D_3(\sigma)$.

Computing the number of permutations with a given $3$-descent set $S$ yields
a few immediate observations. For example, the number of permutations
$\sigma\in\S_n$ with~$D_3(\sigma)$ equal to a fixed subset $S\subseteq [n-2]$
is divisible by $n$. This fact becomes
clear upon the realization that $D_3(\sigma)$ is preserved when the elements
of $\sigma$ are cyclically shifted, so that $1$ becomes $2$, $2$ becomes $3$,
and so on. As a result, it makes sense to focus on the set $\tilde{\S}_n$ of
permutations of $[n]$ with the 
first element equal to $1$. A second, less trivial
observation arising from early calculations is that the number of permutations
in $\tilde{\S}_n$ whose $3$-descent set is empty is the Euler number
$E_{n-1}$. 

This second observation follows from the equidistribution of the statistic
$D_3$ on the set $\tilde{\S}_{n+1}$ with another variation on the descent
set statistic, this time on $\S_n$, which we call the \emph{alternating descent
set} (Theorem \ref{bijection_alt3}).
It is defined as the set of positions $i$ at which the permutation has
an \emph{alternating descent}, which is a regular descent if $i$ is odd or
an ascent if $i$ is even. Thus the alternating descent set $\hat{D}(\sigma)$
of a permutation $\sigma$ is the set of places where $\sigma$ deviates from
the alternating pattern.

Many of the results in this paper 
that were originally motivated by the generalized
descent statistic $d_3(\sigma) = |D_3(\sigma)|$ are actually given in terms
of the alternating descent statistic $\hat{d}(\sigma) = |\hat{D}(\sigma)|$.
We show that the alternating Eulerian polynomials, defined 
as $\hat{A}_n(t) := \sum_{\sigma\in\S_n} t^{\hat{d}(\sigma)+1}$
by analogy with
the classical Eulerian polynomials, have the generating function
$$
\sum_{n\geq 1} \hat{A}_n(t)\cdot {u^n\over n!} = 
{t\left(1-h\bigl({u(t-1)}\bigr)\right) \over h\bigl(u(t-1)\bigr)-t}
$$
where $h(x) = \tan x + \sec x$, so that the difference with the classical
formula (\ref{des_inv_gf}) below
(specialized at $q=1$) is only in that the exponential
function is replaced by tangent plus secant (Theorem \ref{final_F}).

A similar parallel becomes apparent in our consideration of the analog of the
well known identity 
\begin{equation}\label{eulerian_identity}
{A_n(t)\over (1-t)^{n+1}} = \sum_{m\geq 1} m^n t^m
\end{equation}
for $\hat{A}_n(t)$. Given a formal power
series $f(x) = 1 + \sum_{n\geq 1} a_n x^n/n!$, we define the symmetric function
$$
g_{f,n} := \sum_{\gamma\models n} {n\choose \gamma} \cdot
a_{\gamma_1} a_{\gamma_2}\cdots \cdot M_\gamma,
$$
where $\gamma$ runs over all compositions of $n$, and
$$
M_\gamma := \sum_{i_1 < i_2 < \cdots} x_{i_1}^{\gamma_1} x_{i_2}^{\gamma_2}\cdots .
$$
Then
(\ref{eulerian_identity}) can be written as
$$
{A_n(t)\over (1-t)^{n+1}} = \sum_{m\geq 1} g_{\exp,n}(1^m)\cdot t^m,
$$
and we have
$$
{\hat{A}_n(t)\over (1-t)^{n+1}} = \sum_{m\geq 1} g_{\tan+\sec,n}(1^m)\cdot t^m,
$$
where $1^m$ denotes setting the variables $x_1$, $x_2$, \ldots, $x_m$ to $1$
and the remaining variables to $0$ (Proposition \ref{gfn_substitution}).

In Section \ref{cd-stuff} we discuss the generating function
$\hat{\Psi}(\aaa,\bbb)$ for the
number of permutations in $\S_n$ with a given alternating descent set
$S\subseteq [n-1]$, denoted $\hat{\beta}_n(S)$, which is analogous to the
generating polynomial $\Psi_n(\aaa,\bbb)$ for the regular descent set statistic
mentioned earlier.
The polynomial $\Psi_n(\aaa,\bbb)$ can be expressed as the $cd$-index
$\Phi_n(\ccc,\ddd)$ of the Boolean algebra $B_n$, where $\ccc=\aaa+\bbb$
and $\ddd=\aaa\bbb+\bbb\aaa$. We show that $\hat{\Psi}_n$ can also be written
in terms of~$\ccc$ and $\ddd$ as $\hat{\Phi}_n(\ccc,\ddd) =
\Phi_n(\ccc,\ \ccc^2-\ddd)$ (Proposition \ref{hat_Phi}), and that the sum of
absolute values of the coefficients of this $(\ccc,\ddd)$-polynomial, which
is the evaluation $\Phi_n(1,2)$, is the $n$-th term of a
notable combinatorial sequence counting permutations in $\S_n$
with no consecutive descents and no descent at the end
(Theorem \ref{sum_coefficients}).
This sequence has properties relevant to this work; in
particular, the logarithm of the corresponding exponential generating function
is an odd function, which is a crucial property of both $e^x$ and $\tan x
+ \sec x$ that emerges repeatedly in the derivations of the results mentioned
above. We discuss the similarities with Euler numbers and alternating
permutations in Section \ref{SWG_permutations}.

It is natural to wonder if the variations of descents introduced thus far can
be accompanied by corresponding variations of inversions. For alternating
descents it seems reasonable to consider \emph{alternating inversions} defined
in a similar manner as pairs of indices $i<j$ such that either $i$ is odd and
the elements in positions $i$ and $j$ form a regular inversion, or else $i$
is even and these two elements do \emph{not} form a regular inversion. As for
$3$-descents, we define the accompanying \emph{$3$-inversion statistic},
where a $3$-inversion is defined as the number of pairs of indices $(i,j)$
such that $i+1<j$ and the elements in positions $i$, $i+1$, and $j$, taken in
this order, constitute an odd permutation of size $3$. 
Let $\hat{\imath}(\sigma)$ and $i_3(\sigma)$ be the number of alternating
inversions and $3$-inversions of a permutation $\sigma$, respectively.
We find that the joint distribution of the pair $(\hat{d}, \hat{\imath})$ of
statistics on the set $\S_n$ is identical to the distribution of the
pair~$(d_3,i_3)$ of statistics on the set $\tilde{\S}_{n+1}$
(Theorem \ref{main_bivariate_identity}).

Stanley \cite{StanleyBinomialPosets} derived a generating function
for the joint distribution of the classical descent and inversion statistics
on $\S_n$:
\begin{equation}\label{des_inv_gf}
1 + \sum_{n\geq 1} \sum_{\sigma\in\S_n} t^{d(\sigma)} q^{\inv(\sigma)}\cdot
{u^n\over [n]_q!} = {1-t\over \Exp_q\bigl(u(t-1)\bigr) - t},
\end{equation}
where $\Exp_q(x) = \sum_{n\geq 0} q^{n\choose 2} x^n/[n]_q!$,
and $d(\sigma)$ and $\inv(\sigma)$ denote the number
of descents and inversions of $\sigma$, respectively.
(Another good reference on the subject
is a recent paper \cite{ShareshianWachs} of Shareshian and Wachs.)
It would be nice to produce an analog of the generating function
(\ref{des_inv_gf}) for these descent-inversion pairs, but this task appears
to be challenging, and it is not even clear what form such a generating
function should have, as the $q$-factorials in the denominators of
(\ref{des_inv_gf}) are strongly connected to $q$-binomial coefficients,
which have a combinatorial interpretation of the number of inversions in a
permutation obtained by concatenating two increasing runs of fixed size.
Nevertheless the bivariate polynomial $\hat{A}_n(t,q)
:= \sum_{\sigma\in\S_n} t^{\hat{d}(\sigma)} q^{\hat{\imath}(\sigma)}$ seems
to be of interest, and in Section \ref{Euler-q-analog}
we direct our attention to
the $q$-polynomials that result if we set $t=0$. This special case concerns
up-down permutations and, more precisely, their distribution according to the
number of alternating inversions. For down-up permutations this distribution
is essentially the same, the only difference being the order of the
coefficients in the $q$-polynomial, and for our purposes it turns out to be
more convenient to work with down-up permutations, so we use the distribution
of $\hat{\imath}$ on them to define a $q$-analog $\hat{E}_n(q)$
of Euler numbers. The formal definition we give is
$$
\hat{E}_n(q) := q^{-\lfloor n^2/4 \rfloor}
\sum_{\sigma\in\Alt_n} q^{\hat{\imath}(\sigma)},
$$
where $\Alt_n$ is the set of down-up permutations of $[n]$. The polynomial
$\hat{E}_n(q)$ is monic with constant term equal to the Catalan number
$c_{\lfloor n/2 \rfloor}$ (Proposition \ref{hatEq_facts}), which hints at the
possibility to express $\hat{E}_n(q)$ as the sum of $c_{\lfloor n/2 \rfloor}$
``nice'' polynomials with constant term $1$. We discover
such an expression in the form of a $q$-analog of a beautiful identity
that represents $E_n$ as the sum of weighted Dyck paths of
length~$2\lfloor n/2 \rfloor$.
In this identity we imagine Dyck paths as starting
at $(0,0)$ and ending at $(2\lfloor n/2 \rfloor, 0)$. We set the weight of
an up-step to be the level at which that step is situated
(the steps that touch the ``ground'' are at level $1$, the steps above them
at level $2$, and so on) and the weight of a down-step to be either the
level of the step (for even $n$) or one plus the level of the step
(for odd $n$). 
We set the weight of the path to be the product of the weights of all its
steps.
The sum of the weights taken over all $c_{\lfloor n/2\rfloor}$
paths then equals $E_n$, and if we replace the weight of a step with the
$q$-analog of the respective integer, we obtain $\hat{E}_n(q)$
(Theorem \ref{weighted_path_identity}).

The original $q=1$ version of the above identity provides a curious connection
between Catalan and Euler numbers. A notable difference between these
numbers is in the generating functions: one traditionally considers the
ordinary generating function for the former and the exponential one for the
latter. An interesting and hopefully solvable problem is to find a generating
function interpolating between the two, and a potential solution could be to
use the above $q$-analog $\hat{E}_n(q)$ of Euler numbers to write
$$
H(q,x) := \sum_{n\geq 0} \hat{E}_n(q)\cdot {x^n\over [n]_q!},
$$
so that $H(1,x) = \tan x + \sec x$ and
$$
H(0,x) = \sum_{n\geq 0} c_{\lfloor n/2\rfloor} x^n =
{(1+x)\left(1-\sqrt{1-4x^2}\right)\over 2x^2}.
$$

\section*{Acknowledgments}

This paper is part of the author's Ph.D.\ thesis.
I would like to thank Pavlo Pylyavskyy for his ideas and conversations that led to this work.
I am also grateful to Richard Stanley and Alex Postnikov for helpful discussions.

\section{Variations on the descent statistic}\label{descents}


Let $\S_n$ be the
set of permutations of $[n]=\{1,\ldots,n\}$, and let $\tilde{\S}_{n}$
be the set of permutations $\sigma_1\sigma_2\cdots \sigma_{n}$ of 
$[n]$ such that $\sigma_1=1$. For a permutation
$\sigma=\sigma_1\cdots \sigma_n$, define the \emph{descent set}
$D(\sigma)$ of $\sigma$ by
$D(\sigma)=\{i\ |\
\sigma_i>\sigma_j\} \subseteq [n-1]$, and set
$d(\sigma)=|D(\sigma)|$. 

We say that a permutation $\sigma$ has a \emph{$3$-descent} at
position $i$ if the permutation $\sigma_i\sigma_{i+1}\sigma_{i+2}$,
viewed as an element of $\S_3$, is odd. Let $D_3(\sigma)$ be the set
of positions at which a permutation $\sigma$ has a $3$-descent, and
set $d_3(\sigma)=|D_3(\sigma)|$. An important property of the
$3$-descent statistic is the following.

\begin{lemma}\label{cyclic_shift}
Let $\omega^c_n$ be the cyclic permutation $(2\ 3\ \ldots\ n\ 
1)$, and let $\sigma\in\S_n$. Then $D_3(\sigma)=D_3(\sigma\omega^c_n)$.
\end{lemma}

\begin{proof}
Multiplying $\sigma$ on the right by $\omega^c_n$ 
replaces each $\sigma_i<n$ by $\sigma_i+1$, and the element of
$\sigma$ equal to $n$ by $1$. Thus the elements of the triples
$\sigma_i\sigma_{i+1}\sigma_{i+2}$ that do not include $n$ maintain
their relative order under this operation, and in the triples that
include $n$, the relative order of exactly two pairs of elements is
altered.
Thus the $3$-descent set of $\sigma$ is preserved.  
\end{proof}

\begin{corollary}\label{corollary_cyclic_shift}
For all $i,j,k,\ell\in [n]$ and $B\subseteq [n-2]$, 
the number of permutations $\sigma\in\S_n$ with $D_3(\sigma)=B$ and
$\sigma_i=j$ is the same as the number of permutations with
$D_3(\sigma)=B$ and $\sigma_k=\ell$.
\end{corollary}

\begin{proof}
The set $\S_n$ splits into orbits of the form
$\{\sigma,\sigma\omega^c_n,
\sigma(\omega^c_n)^2,\ldots,\sigma(\omega^c_n)^{n-1}\}$, and each such subset
contains exactly one permutation with a $j$ in the $i$-th position for
all
$i,j\in [n]$.
\end{proof}

Next, we define another variation on the descent statistic. We say
that a permutation $\sigma=\sigma_1\cdots\sigma_n$ 
has an \emph{alternating descent} at
position $i$ if either $\sigma_i>\sigma_{i+1}$ and $i$ is odd, or else
if
$\sigma_i<\sigma_{i+1}$ and $i$ is even. Let $\hat{D}(\sigma)$ be the
set of positions at which $\sigma$ has an alternating descent, and set
$\hat{d}(\sigma)=|\hat{D}(\sigma)|$.

Our first result relates the last two statistics by asserting that
the $3$-descent sets of permutations in $\tilde{\S}_{n+1}$ are
equidistributed with the alternating descent sets of permutations in
$\S_n$.

\begin{theorem}\label{bijection_alt3}
Let $B\subseteq [n-1]$. The number of permutations
$\sigma\in\tilde{\S}_{n+1}$ with $D_3(\sigma)=B$ is equal to the number
of permutations $\omega\in\S_n$ with $\hat{D}(\omega)=B$.
\end{theorem}

\begin{proof_}{Proof (by Pavlo Pylyavskyy, private communication).}
We construct a bijection between $\tilde{\S}_{n+1}$ and
$\S_n$ mapping permutations with $3$-descent set $B$ to permutations
with alternating descent set $B$.

Start with a permutation in $\sigma \in
\tilde{\mathfrak{S}}_n$. 
We construct the
corresponding permutation $\omega$ in $\mathfrak{S}_n$ by the following
procedure. Consider $n+1$ points on a circle, and label them with
numbers from $1$ to $n+1$ in the clockwise
direction. 
For convenience, we refer to these points by their labels.
For $1\leq i \leq n$, draw a line segment connecting $\sigma_{i}$
and
$\sigma_{i+1}$.
The segment $\sigma_{i}\sigma_{i+1}$ divides the circle into two arcs.
Define the sequence $C_1$, \ldots, $C_n$, where $C_i$ is one of the two
arcs
between $\sigma_{i}$ and $\sigma_{i+1}$, 
according
to the following rule. 
Choose $C_1$ to be the arc between $\sigma_1$ and $\sigma_2$
corresponding to going from $\sigma_1$ to $\sigma_2$ in the
clockwise direction. For $i>1$,
given the choice of $C_{i-1}$, let $C_i$ be the
arc 
between $\sigma_{i}$ and $\sigma_{i+1}$ that
either {\it {contains}} or {\it {is
contained in}} $C_{i-1}$. The choice of such an arc is always
possible and unique. Let $\ell(i)$ denote how many of the $i$ points
$\sigma_1, \ldots, \sigma_{i}$, including $\sigma_{i}$, 
are contained in $C_i$.

Now, construct the sequence of permutations
$\omega^{(i)}=\omega^{(i)}_1 \ldots \omega^{(i)}_i\in\S_i$,
$1\leq i \leq n$,
as follows. Let $\omega^{(1)} = \ell(1)$. Given $\omega^{(i-1)}$, 
set $\omega^{(i)}_i =
\ell(i)$, and let $\omega^{(i)}_1\ldots \omega^{(i)}_{i-1}$ be the
permutation
obtained from $\omega^{(i-1)}$ by adding $1$ to all elements 
which are greater than or equal to $\ell_i$.
Finally, set $\omega=\omega^{(n)}$.

Next, we argue that the map $\sigma \mapsto \omega$ 
is a bijection. Indeed, from the subword $\omega_1
\ldots \omega_i$ of $\omega$ one can recover $\ell(i)$ since
$\omega_i$ is the $\ell(i)$-th smallest element of the set
$\{\omega_1,\ldots,\omega_i\}$.
Then one can reconstruct one by one the arcs $C_i$ and the
segments connecting $\sigma_{i}$ and $\sigma_{i+1}$ as follows. If
$\ell(i)>\ell(i-1)$ then $C_i$ contains $C_{i-1}$, and if $\ell(i) \leq
\ell(i-1)$ then $C_i$ is contained in $C_{i-1}$. Using this observation and
the number $\ell(i)$ of the points $\sigma_1, \ldots, \sigma_{i}$
contained in $C_i$, one can determine the position of the point $\sigma_{i+1}$
relative to the points $\sigma_1,\ldots, \sigma_{i}$.

It remains to check that $D_3(\sigma)=\hat{D}(\omega)$. Observe that $\sigma$
has a
$3$-descent in position $i$ if and only if 
the triple of points $\sigma_i, \sigma_{i+1}, \sigma_{i+2}$ on the
circle is oriented counterclockwise. Also, observe that
$\omega_i>\omega_{i-1}$ if and only if $C_{i-1}\subset C_i$. Finally,
note that
$C_{i-1} \subset C_i \supset C_{i+1}$ or $C_{i-1} \supset C_i
\subset C_{i+1}$
if and only if triples $\sigma_{i-1}, \sigma_{i}, \sigma_{i+1}$ and 
$\sigma_{i}, \sigma_{i+1}, \sigma_{i+2}$ have
the same orientation.
We now show by induction on $i$ that $i\in D_3(\sigma)$ if and only if
$i\in \hat{D}(\omega)$.
From the choice of $C_1$ and $C_2$, it follows that $C_1\subset C_2$
if and only if $\sigma_3>\sigma_2$, and hence $\omega$ has an
(alternating) descent at position $1$ if and only if
$\sigma_1\sigma_2\sigma_3 = 1\sigma_2\sigma_3$ is an odd permutation.
Suppose the claim holds for $i-1$. By the above observations,
we have $\omega_{i-1} < \omega_i > \omega_{i+1}$ or 
$\omega_{i-1} > \omega_i < \omega_{i+1}$ if and
only if the permutations $\sigma_{i-1}\sigma_{i}\sigma_{i+1}$ and
$\sigma_{i}
\sigma_{i+1} \sigma_{i+2}$ have the same sign. In other words, $i-1$
and
$i$ are either both contained or both not contained in 
$\hat{D}(\omega)$
if and only if they are either both contained or
both not contained in 
$D_3(\sigma)$. It follows that
$i\in D_3(\sigma)$ if and only if $i\in\hat{D}(\omega)$.
\end{proof_}

An important special case of Theorem \ref{bijection_alt3} is
$B=\varnothing$.
A permutation $\sigma\in\S_n$ has $\hat{D}(\sigma)=\varnothing$ if and
only if it is an \emph{alternating (up-down)} permutation,
i.e. $\sigma_1 < \sigma_2 > \sigma_3 < \cdots$. The number of such
permutations of size $n$ is the \emph{Euler number} $E_n$. Thus we get
the following corollary:

\begin{corollary}\label{corollary_Euler}
(a) The number of permutations in $\tilde{\S}_{n+1}$ with no $3$-descents
is $E_n$.

\vskip4pt
\noindent
(b) The number of permutations in $\S_{n+1}$ with no $3$-descents is 
$(n+1)E_n$.
\end{corollary}

\begin{proof}
Part (b) follows from Corollary \ref{corollary_cyclic_shift}: for each
$j\in
[n+1]$, there are $E_n$ permutations in $\S_{n+1}$ beginning with $j$.
\end{proof}

Permutations with no $3$-descents can be equivalently described as
simultaneously avoiding \emph{generalized} patterns $132$, $213$, and $321$
(meaning, in this case, triples of \emph{consecutive} elements
with one of these relative orders). Corollary \ref{corollary_Euler}(b)
appears in the paper \cite{KitaevMansour} of Kitaev and Mansour
on simultaneous avoidance of generalized patterns. Thus the above
construction yields a bijective proof of their result.

\section{Variations on the inversion statistic}\label{inversions}

In this section we introduce analogs of the inversion statistic on permutations
corresponding to the $3$-descent and the alternating descent statistics 
introduced in Section \ref{descents}. First, let us recall the
standard inversion statistic. For $\sigma\in\S_n$, let $a_i$ be the
number of indices $j>i$ such that $\sigma_i > \sigma_j$, and set $\code(\sigma)
= (a_1,\ldots,a_{n-1})$ and $\inv(\sigma) = a_1+\cdots + a_{n-1}$.

For a permutation $\sigma\in\S_n$ and $i\in [n-2]$, 
let $c^3_i(\sigma)$ be the number of indices $j>i+1$ such that
$\sigma_i \sigma_{i+1} \sigma_j$ is an
odd permutation,
and set
$\code_3(\sigma)=(c^3_1(\sigma),c^3_2(\sigma),\ldots,c^3_{n-2}(\sigma))$.
Let $C_k$ be the set of $k$-tuples
$(a_1,\ldots,a_k)$ of non-negative integers
such that $a_i \leq k+1-i$. Clearly,  
$\code_3(\sigma)\in C_{n-2}$. 

\begin{lemma}\label{cyclic_shift_inversions}
Let $\omega^c_n$ be the cyclic permutation $(2\ 3\ \ldots\ n\ 
1)$, and let $\sigma\in\S_n$. Then $\code_3(\sigma)=\code_3(\sigma\omega^c_n)$.
\end{lemma}

\begin{proof}
The proof is analogous to that of Lemma \ref{cyclic_shift}.
\end{proof}

\begin{proposition}\label{code_3}
The restriction $\code_3 : \tilde{\S}_n \rightarrow C_{n-2}$ is a bijection.
\end{proposition}

\begin{proof}
Since $|\tilde{\S}_n|=|C_{n-2}| = (n-1)!$, it suffices to show that
the restriction of $\code_3$ to $\tilde{\S}_n$
is surjective.
We proceed by induction on
$n$. The claim is trivial for $n=3$. Suppose it is true for $n-1$,
and let $(a_1,\ldots,a_{n-2})\in C_{n-2}$. 
Let $\tau$ be the unique permutation in $\tilde{\S}_{n-1}$ such that
$\code_3(\tau)=(a_2,\ldots,a_{n-2})$. 
For $1\leq\ell\leq n$, let $\ell*\tau$ be the
permutation in $\S_n$ beginning with $\ell$ such that the relative order
of last $n-1$ elements of $\ell*\tau$ is the same as that of the elements
of $\tau$. Setting $\ell=n-a_1$ we obtain $\code_3(\ell*\tau) =
(a_1,\ldots,a_{n-2})$ since $\ell\ 1\ m$ is an odd permutation if and
only if $\ell < m$, and there are exactly $a_1$ elements of
$\ell*\tau$ that are greater than $\ell$. Finally, by Lemma
\ref{cyclic_shift_inversions}, the permutation $\sigma =
(\ell*\tau)(\omega^c_n)^{1-a_1}\in \tilde{\S}_{n}$ satisfies
$\code_3(\sigma)=(a_1,\ldots,a_{n-2})$.   
\end{proof}

Let $i_3(\sigma) = c^3_1(\sigma) + c^3_2(\sigma) +
\cdots
+ c^3_{n-2}(\sigma)$. An immediate consequence of Proposition
\ref{code_3} is that $i_3(1*\sigma)$ is 
a \emph{Mahonian statistic} on permutations $\sigma\in\S_{n}$:

\begin{corollary}\label{code3_mahonian}
We have
$$
\sum_{\sigma\in\S_n} q^{i_3(1*\sigma)} = (1+q)(1+q+q^2)\cdots
(1+q+q^2+\cdots+q^{n-1}).
$$
\end{corollary}

For a permutation $\sigma\in\S_n$ and $i\in [n-1]$, define 
$\hat{c}_i(\sigma)$ to be the number of indices $j>i$ such that
$\sigma_i>\sigma_j$ if $i$ is odd, or the number of indices $j>i$ such
that $\sigma_i<\sigma_j$ if $i$ is even.
Set $\hat{\code}(\sigma) =
(\hat{c}_1(\sigma), \ldots, \hat{c}_{n-1}(\sigma))\in C_{n-1}$ and
$\hat{\imath}(\sigma) = \hat{c}_1(\sigma) + \cdots + \hat{c}_{n-1}(\sigma)$.

\begin{proposition}\label{hat_code}
The map $\hat{\code} : \S_n \rightarrow C_{n-1}$ is a bijection.
\end{proposition}

\begin{proof}
The proposition follows easily from the fact that if $\code(\sigma) =
(a_1,\ldots, a_{n-1})$
is the standard inversion code of $\sigma$, then
$\hat{\code}(\sigma)=(a_1,n-2-a_2,a_3,n-4-a_4,\ldots)$. Since the
standard inversion code is a bijection between $\S_n$ and $C_{n-1}$,
so is $\hat{\code}$. 
\end{proof}

\begin{corollary}\label{hat_code_mahonian}
We have
$$
\sum_{\sigma\in\S_n} q^{\hat{\imath}(\sigma)} = (1+q)(1+q+q^2)\cdots
(1+q+q^2+\cdots+q^{n-1}).
$$
\end{corollary}

Another way to deduce Corollary \ref{hat_code_mahonian} is via the
bijection
$\sigma \leftrightarrow \sigma^\vee$, where
$$
\sigma^\vee = \sigma_1 \sigma_3 \sigma_5 \cdots
\sigma_6 \sigma_4 \sigma_2.
$$

\begin{proposition}\label{sigma_check}
We have $\hat{\imath}(\sigma) = \inv(\sigma^\vee)$.
\end{proposition}

\begin{proof}
It is easy to verify that a pair $(\sigma_i, \sigma_j)$, $i<j$, 
contributes to
$\hat{\imath}(\sigma)$ if and only if it contributes to $\inv(\sigma^\vee)$.
\end{proof}

Next, we prove a fundamental relation between the variants of the
descent and the inversion statistics introduced thus far.

\begin{theorem}\label{main_bivariate_identity}
We have
$$
\sum_{\sigma\in \tilde{\S}_{n+1}} t^{d_3(\sigma)} q^{i_3(\sigma)}
= \sum_{\omega\in \S_n} t^{\hat{d}(\omega)} q^{\hat{\imath}(\omega)}.
$$
\end{theorem}

\begin{proof}
The theorem is a direct consequence of the following proposition.

\begin{proposition}\label{equal_codes}
If $\code_3(\sigma) = \hat{\code}(\omega)$ for some
$\sigma\in\tilde{\S}_{n+1}$ and $\omega\in\S_n$, then
$D_3(\sigma)=\hat{D}(\omega)$. 
\end{proposition}

\begin{proof}
The alternating descent set of $\omega$ can be obtained from
$\hat{\code}(\omega)$ as follows:

\begin{lemma}\label{hatD_from_hatcode}
For $\omega\in\S_n$, 
write $(a_1,\ldots,a_{n-1})=\hat{code}(\omega)$, and set $a_n=0$.
Then
$\hat{D}(\omega)=\{i\in [n-1]\ |\ a_i+a_{i+1} \geq n-i\}$.
\end{lemma}

\begin{proof}
Suppose $i$ is odd; then 
if $\omega_i>\omega_{i+1}$, 
i.e. $i\in\hat{D}(\omega)$,
then for each $j>i$ we have
$\omega_i>\omega_j$
or $\omega_{i+1}<\omega_j$ or both, so $a_i+a_{i+1}$ is not smaller
than $n-i$, which is
the
number of elements of $\omega$ to the right of $\omega_i$; if on the
other hand $\omega_i<\omega_{i+1}$, i.e. $i\notin\hat{D}(\omega)$,
then for each $j>i$, at most one of the inequalities
$\omega_i>\omega_j$ and $\omega_{i+1}<\omega_j$ holds, and neither
inequality holds for $j=i+1$, so $a_i+a_{i+1}\leq n-i-1$, which is the
number of elements of $\omega$ to the right of $\omega_{i+1}$. The
case of even $i$ is analogous.
\end{proof}

We now show that the $3$-descent set of $\sigma$ can be obtained from
$(a_1,\ldots,a_{n-1})$ in the same way.

\begin{lemma}\label{D3_from_code3}
For $\sigma\in\tilde{\S}_{n+1}$, 
write $(a_1,\ldots,a_{n-1})=\code_3(\sigma)$, and set $a_n=0$.
Then
$D_3(\sigma)=\{i\in [n-1]\ |\ a_i+a_{i+1} \geq n-i\}$.
\end{lemma}

\begin{proof}
Let $B=D_3(\sigma)$, and let 
$\sigma' = \sigma (\omega_{n+1}^c)^{1-\sigma_i} \in
\S_{n+1}$. Then
$\sigma'_i=1$, and by Lemmas \ref{cyclic_shift} and
\ref{cyclic_shift_inversions}, we have $D_3(\sigma')=D_3(\sigma)=B$ and
$\code_3(\sigma')=\code_3(\sigma)$. 

Suppose that $1=\sigma'_i<\sigma'_{i+1}<\sigma'_{i+2}$. Then $i\notin B$,
and for each $j>i+2$, at most one of the permutations
$\sigma_i'\sigma'_{i+1} \sigma'_{j} = 1\sigma'_{i+1}\sigma'_j$ and
$\sigma'_{i+1}  \sigma'_{i+2}  \sigma'_j$ is odd, because 
$1 \sigma'_{i+1} \sigma'_{j}$ is odd if and only if
$\sigma'_{i+1} > \sigma'_j$, and $\sigma'_{i+1} \sigma'_{i+2}
\sigma'_j$ is odd if and only if $\sigma'_{i+1} < \sigma'_j <
\sigma'_{i+2}$. Hence $a_i+a_{i+1}$ is at most $n-1-i$, which is
the number of indices $j\in [n+1]$ such that $j>i+2$.

Now suppose that $1=\sigma'_i < \sigma'_{i+1} > \sigma'_{i+2}$. Then
$i\in B$, and for each $j>i+2$, at least one of the permutations
$\sigma'_i \sigma'_{i+1} \sigma'_{j} = 1 \sigma'_{i+1} \sigma'_j$ 
and $\sigma'_{i+1}
\sigma'_{i+2}  \sigma'_{j}$ is odd, because $\sigma'_{i+1}
> \sigma'_j$ makes $1 \sigma'_{i+1} \sigma'_{j}$ odd, and
$\sigma'_{i+1} < \sigma_j$ makes $\sigma'_{i+1} \sigma'_{i+2}
\sigma'_j$ odd. Thus each index $j>i+1$ contributes to at least
one of $a_i$ and $a_{i+1}$, so $a_i+a_{i+1} \geq n-i$, which is the
number of indices $j\in [n+1]$ such that $j > i+1$.
\end{proof}

Proposition \ref{equal_codes} follows from Lemmas
\ref{hatD_from_hatcode} and \ref{D3_from_code3}. 
\end{proof}

Combining the results of the above discussion, we conclude that both
polynomials of Theorem \ref{main_bivariate_identity} are equal to
$$
\sum_{(a_1,\ldots,a_{n-1})\in C_{n-1}}
t^{|\hat{D}(a_1,\ldots,a_{n-1})|}\  
q^{a_1+\cdots + a_{n-1}},
$$
where $\hat{D}(a_1,\ldots,a_{n-1}) = \{i\in [n-1]\ |\ a_i+a_{i+1} \geq n-i\}$.
\end{proof}

Note that the bijective correspondence
$$
\sigma\in\S_n
\xrightarrow{\ \ \hat{\code}\ \ }
c\in C_{n-1}
\xrightarrow{\ \ (\code_3)^{-1}\ \ }
\omega\in\tilde{\S}_{n+1}
$$
satisfying $\hat{D}(\sigma) = D_3(\omega)$ yields another bijective
proof of Theorem \ref{bijection_alt3}.


Besides the inversion statistic, the most famous Mahonian statistic on
permutations is
the \emph{major index}. For $\sigma \in \S_n$, define the major index
of $\sigma$ by
$$
\maj(\sigma) = \sum_{i\in D(\sigma)} i.
$$
Our next result reveals a close relation between the major index and
the $3$-inversion statistic $i_3$.

\begin{proposition}\label{i3_bmaj}
For $\sigma\in \S_n$, write $\sigma^r = \sigma'_n\cdots
\sigma'_2\sigma'_1$, where $\sigma'_i = n+1-\sigma_i$. 
Then 
$$i_3(1*\sigma) = \maj(\sigma^r).$$
\end{proposition}

\begin{proof}
Let $\sigma = 1*\omega \in \tilde{\S}_{n+1}$.
Let $D(\sigma) = \{b_1 < \cdots < b_d\}$. Write $\sigma = \tau^{(1)}
\tau^{(2)} \cdots \tau^{(d+1)}$, where 
$\tau^{(k)} = \sigma_{b_{k-1}+1}
\sigma_{b_{k-1}+2} \cdots \sigma_{b_k}$ and $b_0=0$ and $b_{d+1}=n$.
In other words, we split $\sigma$ into ascending runs between
consecutive descents. 
Fix an element $\sigma_j$ of $\sigma$, 
and suppose $\sigma_j \in \tau^{(k)}$. We claim that there are exactly
$k-1$ indices $i < j-1$ such that $\sigma_i \sigma_{i+1} \sigma_j$
is an odd permutation. For each ascending run $\tau^{(\ell)}$,
$\ell<k$, there is at most one element $\sigma_i\in\tau^{(\ell)}$ such
that $\sigma_i < \sigma_j < \sigma_{i+1}$, in which case
$\sigma_i\sigma_{i+1}\sigma_j$ is odd. There is no such element in
$\tau^{(\ell)}$ if and only if the first element
$\sigma_{b_{\ell-1}+1}$ of $\tau^{(\ell)}$ is greater than $\sigma_j$,
or the last element $\sigma_{b_\ell}$ of $\tau^{(\ell)}$ is smaller
than $\sigma_j$. In the former case we have $\sigma_{b_{\ell}-1} >
\sigma_{b_{\ell}} > \sigma_j$, so $\sigma_{b_{\ell}-1}
\sigma_{b_{\ell}} 
\sigma_j$ is odd, and in the latter case, $\sigma_j > \sigma_{b_\ell} >
\sigma_{b_\ell + 1}$, so $\sigma_{b_\ell} 
\sigma_{b_\ell + 1} \sigma_j$ is odd. 
Thus we 
obtain a one-to-one
correspondence between the $k-1$ ascending runs $\tau^{(1)}$, \ldots,
$\tau^{(k-1)}$ and elements $\sigma_i$ such that
$\sigma_i\sigma_{i+1}\sigma_j$ is an odd permutation.

We conclude that for each $\tau^{(k)}$, there are $(k-1)\cdot (b_k -
b_{k-1})$ odd triples $\sigma_i \sigma_{i+1} \sigma_j$ with
$\sigma_j\in\tau^{(k)}$, and hence
$$
i_3(\sigma) = \sum_{k=1}^{d+1} (k-1)\cdot (b_k -
b_{k-1}) = 
$$
$$
= (b_{d+1} - b_d) + (b_{d+1} - b_d + b_d - b_{d-1}) + (b_{d+1} - b_d +
b_d - b_{d-1} + b_{d-1} - b_{d-2}) + \cdots =
$$
$$
=
\sum_{m=1}^d (n-b_m).
$$
We have $D(\omega) = \{b_1-1, b_2-1, \ldots, b_d-1\}$, from where
it is not hard to see that $D(\omega^r) = \{n-b_d, n-b_{d-1}, \ldots,
n-b_1\}$. The proposition follows.
\end{proof}

Observe that for a permutation $\pi$ with
$\pi^r = \pi'_m \cdots \pi'_1$, the triple $\pi_{i}
\pi_{i+1} \pi_{i+2}$ is odd if and only if the triple $\pi_{i+2}
\pi_{i+1} \pi_i$ is even, which in turn is the case if and
only if the triple $\pi'_{i+2}\pi'_{i+1}\pi'_i$
of consecutive elements of $\pi^r$ is odd. Thus
$d_3(\pi) = d_3(\pi^r)$, and we obtain the following corollary.

\begin{corollary}\label{other_bivariate_identity}
We have
$$
\sum_{\sigma\in\tilde{\S}_{n+1}} t^{d_3(\sigma)} q^{i_3(\sigma)}
= \sum_{\omega\in\S_n} t^{d_3(\omega\circ (n+1))} q^{\maj(\omega)},
$$
where $\omega\circ (n+1)$ is the permutation obtained by appending
$(n+1)$ to $\omega$.
\end{corollary}

\begin{proof}
To deduce the identity from Proposition \ref{i3_bmaj}, write $\sigma =
1*\pi$ and set $\omega = \pi^r$, so that $\omega\circ (n+1) = \sigma^r$.
\end{proof}

In the language of permutation patterns, the statistic
$i_3(\sigma)$ can be defined as the total number of occurrences of generalized
patterns $13$-$2$, $21$-$3$, and $32$-$1$ in $\sigma$.
(An occurrence of a generalized pattern $13$-$2$ in a permutation
$\sigma = \sigma_1 \sigma_2 \cdots$ is a pair of indices $(i,j)$ such that
$i +1 < j$ and
$\sigma_i$, $\sigma_{i+1}$, and $\sigma_{j}$ have the same relative order as
$1$, $3$, and $2$, that is, $\sigma_{i} < \sigma_{j} < \sigma_{i+1}$, and the
other two patterns are defined analogously.) In \cite{BabsonSteingrimsson}
Babson and Steingr\'{i}msson mention the Mahonian statistic 
$\STAT(\sigma)$, which is defined as $i_3(\sigma)$ (treated in terms of
the aforementioned patterns) plus $d(\sigma)$. In the permutation
$\sigma\circ (n+1)$, where $\sigma\in\S_n$, the descents of $\sigma$ and the 
last element $n+1$ constitute all occurrences of the pattern $21$-$3$
involving $n+1$, and hence
$i_3\bigl(\sigma\circ (n+1)\bigr) = \STAT(\sigma)$. 

\section{Variations on Eulerian polynomials}\label{Eulerian-polynomials}

Having introduced two new descent statistics, it is natural to look at the
analog of the Eulerian polynomials representing their common distribution on
$\S_n$. First, recall the definition of the classical $n$-th Eulerian 
polynomial:
$$
A_n(t) := \sum_{\sigma\in\S_n} t^{d(\sigma)+1} = \sum_{k=1}^n A(n,k)\cdot t^k,
$$
where $A(n,k)$ is the number of permutations in $\S_n$ with $k-1$ descents.
There is a well-known formula for the exponential generating function
for Eulerian polynomials:
\begin{equation}\label{egf_eulerian}
E(t,u) = \sum_{n\geq 1} A_n(t) \cdot {u^n\over n!}
= {t(1-e^{u(t-1)})\over e^{u(t-1)}-t}.
\end{equation}

In this section we consider analogs of Eulerian numbers and
polynomials for our variations of the descent statistic. Define the
\emph{alternating Eulerian polynomials} $\hat{A}_n(t)$ by
$$
\hat{A}_n(t) := \sum_{\sigma\in\S_n} t^{\hat{d}(\sigma)+1}
= \sum_{k=1}^n \hat{A}(n,k)\cdot t^k,
$$
where $\hat{A}(n,k)$ is the number of permutations in $\S_n$ with
$k-1$ alternating descents. Our next goal is to find an expression for
the exponential generating function
$$
F(t,u) := \sum_{n\geq 1} \hat{A}_n(t)\cdot {u^n\over n!}.
$$

We begin by deducing a formula for the number of permutations in
$\S_n$ with a given alternating descent set. For $S\subseteq [n-1]$,
let $\hat{\beta}_n(S)$ be the number of permutations $\sigma\in\S_n$
with $\hat{D}(\sigma)=S$, and let $\hat{\alpha}_n(S) =
\sum_{T\subseteq S} \hat{\beta}_n(T)$ be the number of
permutations $\sigma\in\S_n$
with $\hat{D}(\sigma) \subseteq S$. For $S = \{s_1 <
\cdots < s_k\} \subseteq [n-1]$, let $\co(S)$ be the composition $(s_1,
s_2 - s_1, s_3-s_2, \ldots, s_k - s_{k-1}, n-s_k)$ of $n$, and for a
composition $\gamma =(\gamma_1, \ldots, \gamma_\ell)$ of $n$, let $S_\gamma$ be
the subset $\{\gamma_1, \gamma_1+\gamma_2, \ldots, \gamma_1 + \cdots +
\gamma_{\ell-1}\}$ of $[n-1]$. Also, define
$$
{n \choose \gamma} := {n\choose \gamma_1,\ldots,\gamma_\ell} = {n!\over
  \gamma_1!\cdots \gamma_\ell!}
$$
and
$$
{n \choose \gamma}_E := {n\choose \gamma_1,\ldots,\gamma_\ell} \cdot
E_{\gamma_1} \cdots E_{\gamma_\ell}.
$$

\begin{lemma}\label{alternating_descent_set}
We have
$$
\hat{\alpha}_n(S) = {n \choose \co(S)}_E
$$
and
$$
\hat{\beta}_n(S) = \sum_{T\subseteq S} (-1)^{|S-T|} {n\choose
  \co(T)}_E.
$$
\end{lemma}

\begin{proof}
Let $S=\{s_1 < \cdots < s_k\}\subseteq [n-1]$. Set $s_0=0$ and $s_{k+1}=n$
for convenience. The alternating descent
set of a permutation $\sigma\in\S_n$ is contained in $S$ if and only
if for all $1 \leq i \leq k+1$, the subword $\tau_i=\sigma_{s_{i-1}+1}
\sigma_{s_{i-1}+2} \cdots \sigma_{s_i}$ forms either an up-down (if
$s_{i-1}$ is even) or a down-up (if $s_{i-1}$ is odd)
permutation. Thus to construct a permutation $\sigma$ with
$\hat{D}(\sigma) \subseteq S$, one must choose one of the ${n\choose
  s_1-s_0, s_2-s_1,\ldots,s_{k+1}-s_{k}} = {n\choose \co(S)}$ ways to
distribute the elements of $[n]$ among the subwords $\tau_1$, \ldots,
$\tau_{k+1}$, and then for each $i\in [k+1]$, choose one of the
$E_{s_i-s_{i-1}}$ ways of ordering the elements within the subword
$\tau_i$. The first equation of the lemma follows. The second equation
is obtained from the first via the inclusion-exclusion principle.
\end{proof}

Now consider the sum
\begin{equation}\label{summation1}
\sum_{S\subseteq [n-1]} {n\choose \co(S)}_E x^{|S|} =
\sum_{S\subseteq [n-1]} \hat{\alpha}_n(S)\cdot x^{|S|} =
\sum_{\sigma\in\S_n} \left(\sum_{T\supseteq \hat{D}(\sigma)} x^{|T|}\right)
\end{equation}
(a permutation $\sigma$ contributes to $\hat{\alpha}_n(T)$ whenever 
$T \supseteq \hat{D}(\sigma)$). The right hand side of
(\ref{summation1}) is equal to 
\begin{eqnarray}
\nonumber
\sum_{\sigma\in\S_n} \sum_{T\supseteq \hat{D}(\sigma)}
x^{\hat{d}(\sigma) + |T-\hat{D}(\sigma)|}
&=& \sum_{\sigma\in\S_n} x^{\hat{d}(\sigma)} 
\sum_{i=0}^{n-1-\hat{d}(\sigma)} {n-1-\hat{d}(\sigma)\choose i}\ x^i \\
&=&
\sum_{\sigma\in\S_n} x^{\hat{d}(\sigma)} (1+x)^{n-1-\hat{d}(\sigma)},
\label{summation2}
\end{eqnarray}
as there are ${n-1-\hat{d}(\sigma)\choose i}$ subsets of $[n-1]$
  containing $\hat{D}(\sigma)$.
Continuing with the right hand side of (\ref{summation2}), we
  get
\begin{equation}\label{summation3}
{(1+x)^n\over x} \cdot\sum_{\sigma\in\S_n} \left({x\over
    1+x}\right)^{\hat{d}(\sigma)+1}
= {(1+x)^n\over x} \cdot\hat{A}_n \left({x\over 1+x}\right).
\end{equation}
Combining equations (\ref{summation1})--(\ref{summation3}), we obtain
\begin{equation}\label{summation4}
\sum_{n\geq 1} \left(\sum_{S\subseteq [n-1]} {n\choose \co(S)}_E
  x^{|S|}\right)\cdot
{y^n\over n!} = 
{1\over x}\cdot
\sum_{n\geq 1} \hat{A}_n \left({x\over
    1+x}\right) \cdot {y^n(1+x)^n\over n!} .
\end{equation}
Since $S \mapsto \co(S)$ is a bijection between $[n-1]$ and the set of
compositions of $n$, 
the left hand side of (\ref{summation4}) is
\begin{equation}\label{summation5}
\sum_{n\geq 1} 
\left(
\sum_{\gamma} {E_{\gamma_1} \cdots E_{\gamma_\ell}\over
  \gamma_1! \cdots \gamma_\ell!}\cdot x^{\ell-1}
\right)\cdot y^n
= {1\over x} \cdot
\sum_{\ell\geq 1} x^\ell\cdot
\left(\sum_{i\geq 1} {E_i y^i \over i!}\right)^\ell,
\end{equation}
where the inside summation in the left hand side is over all
compositions $\gamma=(\gamma_1,\ldots,\gamma_\ell)$ of $n$. Applying
the well-known formula $\sum_{j\geq 0} E_j y^j/ j! = \tan y + \sec y$,
the right hand side of (\ref{summation5}) becomes
\begin{equation}\label{summation6}
{1\over x}\cdot \sum_{\ell\geq 1} x^\ell (\tan y + \sec y - 1)^\ell
= {1\over x} \cdot \left({1\over 1 - x (\tan y + \sec y - 1)} -1 \right).
\end{equation}
Now set $t={x\over 1+x}$ and $u = y(1+x)$. Equating the right hand sides of
(\ref{summation4}) and (\ref{summation6}), we obtain
\begin{equation}\label{summation7}
F(t,u) = \sum_{n\geq 1} \hat{A}_n(t)\cdot {u^n\over n!} = 
\left({1\over 1 - x (\tan y + \sec y - 1)} -1 \right). 
\end{equation}
Finally, applying the inverse substitution $x = {t\over 1-t}$ and $y=u(1-t)$
and simplifying
yields
an expression for $F(t,u)$:
\begin{eqnarray}
\nonumber
F(t,u)&=&{x(\tan y + \sec y -1)
\over 1 - x(\tan y + \sec y -1)}\\
\nonumber
&=&
{t\over 1-t}\cdot \left({\tan y + \sec y - 1
\over 1 - {t\over 1-t} \cdot (\tan y + \sec y - 1)}\right) \\
&=& {t\cdot \bigl(\tan (u(1-t)) + \sec (u(1-t)) - 1\bigl) \over
1 - t\cdot\bigl(\tan (u(1-t)) + \sec (u(1-t))\bigl)}.
\label{almost_final_F}
\end{eqnarray}

Using the property $(\tan z + \sec z)(\tan (-z) + \sec (-z))=1$, we can 
rewrite the above expression for $F(t,u)$ as follows:

\begin{theorem}\label{final_F}
We have
$$
F(t,u) = {t\cdot\bigl( 1 - \tan(u(t-1)) - \sec(u(t-1))\bigl)
\over \tan(u(t-1)) + \sec(u(t-1)) - t}.
$$
\end{theorem}
Thus $F(t,u)$ can be expressed by replacing the exponential function in the
formula~(\ref{egf_eulerian}) for $E(t,u)$ by tangent plus secant. In fact,
omitting the Euler numbers and working with standard multinomial coefficients
gives a proof of (\ref{egf_eulerian}).

A basic result on Eulerian polynomials is the identity
\begin{equation}\label{A_over_identity}
{A_n(t)\over (1-t)^{n+1}} = \sum_{m\geq 1} m^n t^m.
\end{equation}
Our next result is a similar identity involving alternating Eulerian 
polynomials. For a partition $\lambda$ of $n$ with $r_i$ parts equal
to $i$, define
$$
z_\lambda := 1^{r_1}\cdot r_1! \cdot 2^{r_2}\cdot r_2! \cdot \cdots.
$$

\begin{theorem}\label{f-eulerian}
Let
$$
\hat{f}_n(m) = \sum_{\lambda} {n!\over z_\lambda}\cdot
{E_{\lambda_1-1} E_{\lambda_2-1} \cdots \over (\lambda_1-1)!
  (\lambda_2-1)! \cdots} \cdot m^{\ell(\lambda)},
$$
where the sum is over all partitions
$\lambda=(\lambda_1,\lambda_2,\ldots,\lambda_{\ell(\lambda)})$ 
of $n$ into odd parts.
Then
$$
{\hat{A}_n(t)\over (1-t)^{n+1}} = \sum_{m\geq 1} \hat{f}_n(m) t^m.
$$
\end{theorem}

\begin{proof}
Let us consider the generating function
$$
G(t,u) := \sum_{n\geq 1} {\hat{A}(t)\over (1-t)^{n+1}}\cdot {u^n\over n!}.
$$
Then, by (\ref{almost_final_F}), we have
\begin{equation}\label{G_expression}
G(t,u) = {1\over 1-t} \cdot F\left(t,\ {u\over 1-t}\right)
= {t\cdot(\tan u + \sec u - 1)\over(1-t)\bigl(1-t\cdot(\tan u + \sec u)\bigr)}.
\end{equation}
Define
$$
H(m,u) := \sum_{n\geq 1} \hat{f}_n(m) \cdot {u^n\over n!}.
$$
This series can be rewritten as follows:
\begin{equation}\label{H_expression_1}
H(m,u) = 
\sum_{n\geq 1} {\hat{f}_n(m)\over n!} \cdot u^n = 
-1+\prod_{i\geq 0}\left(\sum_{j\geq 0} 
{\left({E_{2i} mu^{2i+1}\over (2i+1)!}\right)^j\over j!}
\right).
\end{equation}
Indeed,
for each $i$, the index $j$ in the summation is the number of parts equal
to $2i+1$ in a partition of $n$ into odd parts, and 
it is not hard to check that
the contribution of $j$ parts equal to $2i+1$ to the appropriate terms of
$\hat{f}_n(m)/n!$ is given by the expression inside the summation on the
right. We subtract $1$ to cancel out the empty partition of $0$ counted by
the product on the right but not by $H(m,u)$.
Continuing with the right hand side of (\ref{H_expression_1}), we get
\begin{eqnarray}
\nonumber
H(m,u)+1
&=& \prod_{i\geq 0} \exp\left({E_{2i} mu^{2i+1}\over (2i+1)!}\right)\\
&=& \exp\left(m\sum_{i\geq 0}\left({E_{2i}u^{2i+1}\over (2i+1)!}\right)\right).
\label{H_expression_2}
\end{eqnarray}
The sum appearing in the right hand side of (\ref{H_expression_2}) is the
antiderivative of $\sec u = \sum_{i\geq 0} E_{2i} u^{2i} / (2i)!$ that vanishes
at $u=0$; this antiderivative is $\ln(\tan u + \sec u)$. Therefore
$$
H(m,u) +1= (\tan u + \sec u)^m.
$$
Hence we have
\begin{equation}\label{H_expression_3}
\sum_{m\geq 1} H(m,u)\cdot t^m = 
{(\tan u + \sec u)\cdot t \over 1 - (\tan u + \sec u)\cdot t}
- {1\over 1-t}.
\end{equation}
It is straightforward to verify that the right hand sides of
(\ref{G_expression}) and (\ref{H_expression_3}) agree, and thus
\begin{equation}\label{final_G_H}
\sum_{n\geq 1} {\hat{A}_n(t)\over (1-t)^{n+1}}\cdot {u^n\over n!}
= G(t,u) = \sum_{m\geq 1} H(m,u) t^m
= \sum_{m,n\geq 1} \hat{f}_n(m) t^m\cdot {u^n\over n!}.
\end{equation}
Equating the coefficients of $u^n/n!$ on both sides of (\ref{final_G_H})
completes the proof of the theorem.
\end{proof}

In the terminology of \cite[Sec.\ 4.5]{StanleyEC1}, Theorem \ref{f-eulerian} 
states that the polynomials $\hat{A}_n(t)$ are the
\emph{$\hat{f}_n$-Eulerian polynomials}.

\section{Eulerian polynomials and symmetric functions}
\label{Eulerian_sym}

The results of the previous section can be tied to the theory of symmetric
functions. Let us recall some basics.
For a composition $\gamma = (\gamma_1, \gamma_2, \ldots, \gamma_k)$,
the \emph{monomial quasisymmetric function}
$M_\gamma(x_1, x_2, \ldots)$ is defined by
$$
M_\gamma := \sum_{1\leq i_1 < \cdots < i_k} x_{i_1}^{\gamma_1}
x_{i_2}^{\gamma_2} \cdots
x_{i_k}^{\gamma_k}.
$$
Let $\pi(\gamma)$ denote the partition obtained by rearranging the parts of
$\gamma$ in non-increasing order. Then for a partition $\lambda$,
the \emph{monomial symmetric function}
$m_\lambda(x_1, x_2, \ldots)$ is defined as
$$
m_\lambda := \sum_{\gamma\ :\ \pi(\gamma)=\lambda} M_\gamma.
$$

Let $f(x)$ be a function given by the formal power series
$$
f(x) = 1 + \sum_{n\geq 1} {a_n x^n\over n!}.
$$
Define the symmetric function $g_{f,n}(x_1,x_2,\ldots)$ by
$$
g_{f,n} := \sum_{\gamma \models n} {n\choose \gamma} \cdot a_{\gamma_1}
a_{\gamma_2} \cdots \cdot M_\gamma = \sum_{\lambda \vdash n} 
{n\choose \lambda} \cdot a_{\lambda_1}
a_{\lambda_2} \cdots \cdot m_\lambda,
$$
where by $\gamma\models n$ and $\lambda\vdash n$ we mean that $\gamma$
and $\lambda$ are a composition and a partition of $n$, respectively.
This function can be thought of as the generating function for numbers
like $\alpha_n(S)$ or $\hat{\alpha}_n(S)$ (the number of permutations
$\sigma\in\S_n$ with $D(\sigma)\subseteq S$ or $\hat{D}(\sigma)\subseteq S$,
respectively).
Our first step is to express $g_{f,n}$ in terms of the \emph{power sum
symmetric functions} $p_k(x_1,x_2,\ldots) = \sum x_i^k$.

Consider the generating function
\begin{equation}\label{Gf_definition}
G_f(x_1,x_2,\ldots; u) := \sum_{n\geq 0} g_{f,n}\cdot {u^n\over n!}.
\end{equation}
Then we have
\begin{equation}\label{Gf_expression_1}
G_f = 
\sum_{n\geq 0}
\sum_{\gamma\models n} {a_{\gamma_1} a_{\gamma_2} \cdots\over
\gamma_1! \gamma_2! \cdots}\cdot M_\gamma u^n
= \prod_{i\geq 1} f(x_i u).
\end{equation}
Now let us write
\begin{equation}\label{log_f}
\ln(f(x)) = \sum_{n\geq 1} {b_n x^n\over n!}.
\end{equation}
Then from (\ref{Gf_expression_1}) we have
\begin{equation}\label{Gf_expression_2}
\ln G_f = \sum_{i\geq 1} \ln(f(x_i u))
= \sum_{n\geq 1} b_n p_n(x_1, x_2, \ldots)\cdot {u^n\over n!}.
\end{equation}
Since the power sum symmetric functions $p_\lambda = p_{\lambda_1}
p_{\lambda_2}\cdots$, with $\lambda$ ranging over all partitions of positive
integers, form a basis for the ring of symmetric functions, the transformation
$p_n \mapsto b_n p_n u^n/(n-1)!$, where $u$ is regarded as a scalar, extends to
a homomorphism of this ring. Applying this homomorphism to the well-known
identity
$$
\exp \sum_{n\geq 1} {1\over n}\cdot p_n =
\sum_{\lambda}z_\lambda^{-1} 
p_\lambda,
$$
where $\lambda$ ranges over all partitions of positive integers,
we obtain from (\ref{Gf_expression_2}) that
\begin{eqnarray}
\nonumber
G_f&=&
\exp \sum_{n\geq 1} {1\over n}\cdot\left({b_n p_n u^n\over (n-1)!}\right)\\
&=&\sum_{\lambda}z_{\lambda}^{-1} \cdot {b_{\lambda_1} b_{\lambda_2}\cdots\over
(\lambda_1-1)!(\lambda_2-1)!\cdots} \cdot p_\lambda u^{|\lambda|}.
\label{Gf_final}
\end{eqnarray}
Comparing the coefficients of $u^n$ in (\ref{Gf_definition}) and
(\ref{Gf_final}), we conclude the following:

\begin{proposition}\label{gfn_expression}
For a function $f(x)$ with $f(0) = 1$ and $\ln(f(x)) = \sum_{n\geq 1}
b_n x^n / n!$ we have
$$
g_{f,n} = \sum_{\lambda\vdash n} {n!\over z_\lambda} \cdot
{b_{\lambda_1}b_{\lambda_2}\cdots\over (\lambda_1 - 1)! (\lambda_2 - 1)!\cdots}
\cdot p_\lambda.
$$
\end{proposition}

Two special cases related to earlier discussion are $f(x) = e^x$ and
$f(x) = \tan x + \sec x$. For $f(x) = e^x$, we have $b_1 = 1$, $b_2 = b_3 = 
\cdots = 0$, and hence $g_{f,n} = p_1^n$. In the case of $f(x)= \tan x+\sec x$,
we have
$$
b_i = \left\{
\begin{array}{ll}
E_{i-1} & \mbox{if $i$ is odd}, \\
0       & \mbox{if $i$ is even},
\end{array}
\right.
$$
thus the coefficient at $p_\lambda$ in the expression of Proposition
\ref{gfn_expression} coincides with the coefficient in the term for $\lambda$
in the definition of the polynomial $\hat{f}_n(m)$ of Theorem~\ref{f-eulerian}.
These observations lead to the following restatements of
the classical identity~(\ref{A_over_identity}) and Theorem~\ref{f-eulerian}.

\begin{proposition}\label{gfn_substitution}
Let $g(1^m)$ denote the evaluation of $g(x_1,x_2,\ldots)$ at
$x_1 = x_2 = \cdots = x_m = 1$, $x_{m+1} = x_{m+2} = \cdots = 0$. Then
$$
{A_n(t)\over (1-t)^{n+1}} = \sum_{m\geq 1} g_{\exp,n}(1^m) \cdot t^m
$$
and
$$
{\hat{A}_n(t)\over(1-t)^{n+1}}=
\sum_{m\geq 1}g_{\tan+\sec,n}(1^m)\cdot t^m.
$$
\end{proposition}

\begin{proof}
We have $p_i(1^m) = m$, and hence $p_\lambda(1^m) =
m^{\ell(\lambda)}$.
\end{proof}

It is an interesting problem to prove Proposition \ref{gfn_substitution}
without referring to the results of Section \ref{Eulerian-polynomials}.
Observe that for $\gamma = (\gamma_1, \gamma_2, \ldots, \gamma_k) \models n$,
we have $M_\gamma(1^m) = {m \choose k}$, 
the number of monomials $x_{i_1}^{\gamma_1} 
x_{i_2}^{\gamma_2} \cdots x_{i_k}^{\gamma_k}$ where $1 \leq i_1 < \cdots < i_k
\leq m$, which are the monomials in the definition of $M_\gamma$ that evaluate
to $1$.

It would also be of interest to relate the observations of this section
to Schur functions. One possibility is to consider the following 
generalization of
the \emph{complete homogeneous symmetric function}.
Let $\varphi_f$ be the homomorphism of the ring of symmetric functions 
defined by $p_n \mapsto b_n p_n / (n-1)!$,
where the $b_i$'s are as in equation (\ref{log_f}). Let
$$
h_{f,n} := \sum_{\lambda \vdash n} z_\lambda^{-1} \varphi_f(p_\lambda).
$$
For $f(x) = (1-x)^{-1}$, the homomorphism $\varphi$ is identity,
and $h_{f,n}$ is the standard complete homogeneous symmetric function $h_n$,
defined to be the sum of all monomials in $x_1$, $x_2$, \ldots, of degree $n$.
Then (\ref{Gf_final}) becomes
$$
G_f = \sum_{n\geq 1} h_{f,n} u^n
$$
(we do not really need $u$ here because of homegeneity). We can define
the generalized Schur function $s_{f,\lambda}$, where
$\lambda = (\lambda_1,\lambda_2,\ldots) \vdash n$, by the
\emph{Jacobi-Trudi identity}
$$
s_{f,\lambda} :=\det\Bigl[\ h_{f,\ \lambda_i -i+j}\ \Bigr]_{1\leq i,j\leq n}\ ,
$$
where $h_{f,0} = 1$ and $h_{f,k} = 0$ for $k<0$ (see
\cite[Sec.\ 7.16]{StanleyEC2}). 
What can be said about $s_{f,\lambda}$
for $f(x)=e^x$ and $f(x) = \tan x + \sec x$?

\section{The alternating Eulerian numbers}\label{alt_eulerian_numbers}

In this section we give a recurrence relation that allows to construct a
triangle of alternating Eulerian numbers $\hat{A}(n,k)$ introduced in Section
\ref{Eulerian-polynomials}.
(Recall that $\hat{A}(n,k)$ denotes the number of permutations
in $\S_n$ with $k-1$ alternating descents.)
The first few rows of this
triangle are given in 
Table~2.1.
\begin{table}[h!]\label{alt_eulerian_triangle}
\begin{center}
\begin{tabular}{ccccccccccccccc}
&&&&&&1&&&&&&\\[6pt]
&&&&&1&&1&&&&&\\[6pt]
&&&&2&&2&&2&&&&\\[6pt]
&&&5&&7&&7&&5&&&\\[6pt]
&&16&&26&&36&&26&&16&&\\[6pt]
&61&&117&&182&&182&&117&&61&\\[6pt]
272&&594&&1056&&1196&&1056&&594&&272
\end{tabular}
\end{center}
\caption{Triangle of alternating Eulerian numbers}
\end{table}

The following lemma provides a way to compute alternating Eulerian numbers
given the initial condition $\hat{A}(n,1) = E_n$.

\begin{lemma}\label{alt_Eulerian_recursion}
For $n \geq k\geq 0$ we have
\begin{eqnarray}
\nonumber
&&
\sum_{i=0}^n\sum_{j=0}^k{n\choose i}\cdot\hat{A}(i,j+1)
\cdot\hat{A}(n-i,k-j+1)\\
&&
= (n+1-k)\hat{A}(n, k+1) + (k+1) \hat{A}(n, k+2).
\label{hatA_recursion}
\end{eqnarray}
\end{lemma}

\begin{proof}
First, suppose that $k$ is even.
The left hand side of the equation counts the number of ways to split the
elements of $[n]$ into two groups of sizes $i$ and $n-i$,
arrange the elements in the first and the
second group so
that the resulting permutations have $j$ and $k-j$ alternating descents,
respectively, and writing down the second permutation after the first to form
a permtutation of $[n]$. This permutation has either $k$ or $k+1$ alternating
descents, depending on whether an alternating descent is produced at position
$i$. For a permutation $\sigma\in\S_n$ with $\hat{\imath}(\sigma)=k$,
there are exactly $n+1-k$ ways to produce $\sigma$ by means of the above
procedure, one for every choice of $i \in \hat{D}(\sigma)\cup\{0,n\}$.
Similarly, for $\sigma\in\S_n$ such that $\hat{\imath}(\sigma) = k+1$, there
are exactly $k+1$ ways to produce $\sigma$, one for every choice of
$i\in\hat{D}(\sigma)$. The identity follows.

As for odd $k$, the same argument is valid, except that the quantity
$\hat{A}(n-i,k-j+1)$ in the left hand side should be interpreted as the number
of ways to arrange the elements of the second group to form a permutation with
$k-j$ alternating \emph{ascents}, which become alternating descents when the
two permutations are concatenated. 
\end{proof}

Recall the generating function
$$
F(t,u) = \sum_{n,k\geq 1} \hat{A}(n,k)\cdot {t^k u^n\over n!}
$$
introduced in Section \ref{Eulerian-polynomials}. An alternative way to
express $F(t,u)$ and obtain the result of Theorem \ref{final_F} is by solving a
partial differential equation arising from the recurrence of Lemma
\ref{alt_Eulerian_recursion}.

\begin{proposition}\label{hatA_PDE}
The function $F(t,u)$ is the solution of the partial differential equation
\begin{equation}\label{PDE}
F^2-F=u\cdot {\partial F\over\partial u}
      + (1-t)\cdot{\partial F\over \partial t}
\end{equation}
with the initial condition $F(0,u) = \tan u + \sec u$.
\end{proposition}

\begin{proof}
Since $\hat{A}(n,0)=0$ for all $n$, the left hand side of 
(\ref{hatA_recursion})
is $n!$ times the coefficient of $t^k u^n$ in $\bigl(F(t,u)\bigr)^2$,
which we denote by $[t^k u^n]F^2$. 
The right hand side of (\ref{hatA_recursion}) is
\begin{eqnarray*}
&&
n!\cdot\left(
{\hat{A}(n,k+1)\over (n-1)!}+
{\hat{A}(n,k+1)\over n!}
-{k\hat{A}(n,k+1)\over n!}+
{(k+1)\hat{A}(n,k+2)\over n!}
\right)\\
&&=n!\cdot \Bigl(
[t^k u^n] F_u + [t^k u^n] F - [t^{k-1} u^n] F_t + [t^k u^n] F_t
\Bigr)\\
&&=n!\cdot [t^k u^n]\left(
u F_u + F - t F_t + F_t
\right),
\end{eqnarray*}
where $F_t$ and $F_u$ denote partial derivatives of $F$ with respect to
$t$ and $u$. Equating the above with $n!\cdot [t^k u^n] F^2$ proves
(\ref{PDE}).
\end{proof}

\section{The generating function for the alternating descent set statistic}
\label{cd-stuff}

Besides the generating polynomials for the alternating descent
statistic, another natural generating function to consider is one counting
permutations by their alternating descent set. We begin by stating some
well-known facts about the analogous generating function for the classical
descent set statistic.

Fix a positive integer $n$. For a subset $S\subseteq [n-1]$, define the
monomial $u_S$ in two non-commuting variables $\aaa$ and $\bbb$ by
$u_S = u_1 u_2 \cdots u_{n-1}$, where
$$
u_i = \left\{\begin{array}{ll}\aaa & \mbox{if $i\notin S$,}\\
                              \bbb & \mbox{if $i\in S$.}\end{array}\right.
$$
Consider the generating function
$$
\Psi_n(\aaa,\bbb) := \sum_{S\subseteq [n-1]} \beta_n(S) u_S,
$$
where $\beta_n(S)$ is the number of permutations in $\S_n$ with descent set
$S$. The polynomial~$\Psi_n(\aaa,\bbb)$ is known as the 
\emph{$ab$-index of the Boolean algebra $B_n$}. A remarkable property of
$\Psi_n(\aaa,\bbb)$ (and also of $ab$-indices of a wide class of posets,
including face lattices of polytopes) is that it can be expressed in terms
of the variables $\ccc = \aaa + \bbb$ and $\ddd = \aaa \bbb + \bbb \aaa$.
The polynomial $\Phi_n(\ccc,\ddd)$ defined by $\Psi_n(\aaa,\bbb)
= \Phi_n(\aaa+\bbb,\ \aaa\bbb+\bbb\aaa)$ is called the \emph{$cd$-index of
$B_n$}.

The polynomial $\Phi_n(\ccc,\ddd)$ has positive integer coefficients, for which
several combinatorial interpretations have been found. Here we give one that
will help establish a connection with the alternating descent set statistic.
We proceed with a definition.

\begin{definition}\label{simsun}
A permutation is \emph{simsun} if, for all $k\geq 0$, removing $k$ largest
elements from it results in a permutation with no consecutive descents.
\end{definition}

Let $\SSS_n$ be the set of simsun permutations in $\S_n$ whose last element
is $n$. (Thus $\SSS_n$ is essentially the set of simsun permutations
of $[n-1]$ with an $n$ attached at the end.) It is known that
$|\SSS_n| = E_n$.

For a permutation $\sigma\in\SSS_n$, define the
$(\ccc,\ddd)$-monomial $\cd(\sigma)$ as follows: write out the descent set
of $\sigma$ as a string of pluses and minuses denoting ascents and descents,
respectively, and then replace each occurrence of `` -- + '' by $\ddd$, and each
remaining plus by $\ccc$. This definition is valid because a simsun permutation
has no consecutive descents. For example, consider the permutation
$423516 \in \SSS_6$. Its descent set in the above notation is
`` -- + + -- + '', and thus $\cd(423516) = \ddd\ccc\ddd$.

The simsun permutations provide a combinatorial expression for the $cd$-index
of~$B_n$:
\begin{equation}\label{Phi_via_simsun}
\Phi_n(\ccc,\ddd) = \sum_{\sigma\in\SSS_n} \cd(\sigma).
\end{equation}

Now let us
define the analog of $\Psi_n(\aaa,\bbb)$ for the alternating descent set
statistic:
$$
\hat{\Psi}_n(\aaa,\bbb):= \sum_{S\subseteq [n-1]} \hat{\beta}_n(S) u_S.
$$

\begin{proposition}\label{hat_Phi}
There exists a polynomial $\hat{\Phi}_n(\ccc,\ddd)$ such that
$$\hat{\Phi}_n(\aaa+\bbb,\ \aaa\bbb+\bbb\aaa)=\hat{\Psi}_n(\aaa,\bbb),$$
namely, $\hat{\Phi}_n(\ccc,\ddd) = \Phi_n(\ccc,\ \ccc^2-\ddd)$.
\end{proposition}

\begin{proof}
Note that $\hat{\Psi}_n(\aaa,\bbb)$ is the polynomial obtained from
$\Psi(\aaa,\bbb)$ by switching the letters at even positions in all the
$(\aaa,\bbb)$-monomials. For example, we have $\Psi_3(\aaa,\bbb) = \aaa\aaa
+ 2\aaa\bbb + 2\bbb\aaa + \bbb\bbb$, so $\hat{\Psi}_3(\aaa,\bbb)
= \aaa\bbb + 2\aaa\aaa + 2\bbb\bbb + \bbb\aaa$. In terms of the variables
$\ccc$ and $\ddd$, this operation corresponds to replacing $\ddd
= \aaa\bbb+\bbb\aaa$ with $\aaa\aaa+\bbb\bbb = \ccc^2 - \ddd$, and
$\ccc = \aaa+\bbb$ with either $\aaa+\bbb$ or $\bbb+\aaa$, which in any case
is still equal to $\ccc$.
\end{proof}

The polynomial $\hat{\Phi}_n(\ccc,\ddd)$ has both positive and negative
coefficients, but the
polynomial~$\hat{\Phi}_n(\ccc,\ -\ddd) = \Phi_n(\ccc,\ \ccc^2+\ddd)$
has only positive coefficients. It would be nice to give a combinatorial
interpretation for these coefficients similar to that of the coefficients
of $\Phi_n(\ccc,\ddd)$, so that the coefficients of
$\hat{\Phi}_n(\ccc,\ -\ddd)$
enumerate permutations of a certain kind according to some statistic.
In what follows we show that the sum of the coefficients of
$\hat{\Phi}_n(\ccc,\ -\ddd)$ is equal to the number of permutations containing
no consecutive descents and not ending with a descent. Let $\R_n$ denote the
set of such permutations of $[n]$.

\begin{figure}[h!]                                                         
\begin{center}
\input{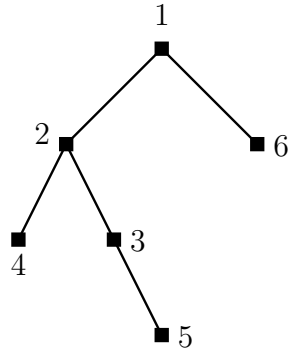}
\end{center}
\caption{The tree $T(423516)$}\label{sample_tree}
\end{figure}

In working with the different kinds of permutations that have emerged thus far
we use the approach of min-tree representation of permutations introduced by
Hetyei and Reiner \cite{HetyeiReiner}.
To a word $w$ whose letters are distinct elements of
$[n]$, associate a labeled rooted planar binary tree according to the following
recursive rule. Let $m$ be the smallest letter of $w$, and write
$w = w_1\circ m\circ w_2$, where $\circ$ denotes concatenation. 
Then form the tree $T(w)$ by labeling the root with $m$ and
setting the left and the right subtrees of the root to be $T(w_1)$ and
$T(w_2)$, respectively. To the empty word we associate the empty tree.
Thus $T(w)$ is an increasing rooted planar binary tree, i.e. the distinction
between left and right children is being made. For example, $T(423516)$ is
the tree shown in Figure \ref{sample_tree}.

To get the word $w$ back from the tree $T(w)$, simply read the labels of the
nodes of $T(w)$ in topological order.

Next, we formulate some of the permutation properties from the above discussion
in terms of the min-tree representation.

\begin{lemma}\label{no_consecutive_descents}
A permutation $\sigma$ has no consecutive descents if and only if the tree
$T(\sigma)$ has no node whose only child is a left child, except maybe for the
rightmost node in topological order.
\end{lemma}

\begin{proof}
Write $\sigma = s_1 s_2 \cdots s_n$ and $T = T(\sigma)$. For convenience,
we refer to the nodes of~$T$ by their labels. We have $s_i > s_{i+1}$ if and
only if $s_{i+1}$ is an ancestor of $s_i$ in~$T$. Since $s_i$ and $s_{i+1}$
are consecutive nodes in the topological reading of $T$, it follows that
$s_{i+1}$ is an ancestor of $s_i$ if and only if $s_i$ has no right child.
Thus we have $s_i > s_{i+1} > s_{i+2}$ if and only if $s_{i+1}$ has no right
child and $s_i$ is a descendant of $s_{i+1}$, i.e. $s_{i+1}$ has a lone left
child. The proposition follows.
\end{proof}

\begin{proposition}\label{Rn_trees}
A permutation $\sigma$ is in $\R_n$ if and only if the tree $T(\sigma)$ has no
node whose only child is a left child.
\end{proposition}

\begin{proof}
We have $s_{n-1} > s_n$ if and only if the rightmost node $s_n$ has a (lone)
left child. The proposition now follows from Lemma
\ref{no_consecutive_descents}.
\end{proof}

\begin{proposition}\label{simsun-trees}
A permutation $\sigma$ is in $\SSS_n$ if and only if the rightmost node
of~$T(\sigma)$ is labeled $n$,
no node has a lone left child, and for every node
$s$ not on the rightmost path (the path from the root to the rightmost node)
that has both a left child $t$ and a right child $u$, the inequality $t>u$
holds.
\end{proposition}

\begin{proof}
If $T(\sigma)$ has a node $s$ not on the rightmost path whose left child $t$
is smaller than its right child $u$, then removing the elements of $\sigma$
that are greater than or equal to $u$ results in a permutation $\sigma'$ such
that in $T(\sigma')$, the node $s$ has a lone left child $t$ and is not the
rightmost node, meaning that $\sigma'$ contains a pair of consecutve descents,
by Lemma \ref{no_consecutive_descents}. If on the other hand $T(\sigma)$ has
no such node $s$, the removing $k$ largest elements of $\sigma$ does not create
any nodes with a lone left child except maybe for the rightmost node.
\end{proof}

One can see that for $\sigma = 423516$, the tree $T(\sigma)$ shown
in Figure~\ref{sample_tree}
satisfies all conditions of Proposition \ref{simsun-trees}, and hence
$423516\in\SSS_6$. Next, we consider the sum of coefficients of
$\hat{\Phi}_n(\ccc,\ -\ddd)$.

\begin{theorem}\label{sum_coefficients}
The sum of coefficients of $\hat{\Phi}_n(\ccc,\ -\ddd)$ is $|\R_n|$.
\end{theorem}

\begin{proof}
The sum of coefficients of $\hat{\Phi}_n(\ccc,\ -\ddd)$ is
$\hat{\Phi}_n(1,-1) = \Phi_n(1,2)$, which equals
$$
\sum_{\sigma\in\SSS_n} 2^{d(\sigma)},
$$
where $d(\sigma)$ is the number of $\ddd$'s in $\cd(\sigma)$, or, equivalently,
the number of descents of $\sigma$. Since the descents of $\sigma$ correspond
to nodes of $T(\sigma)$ that have no right child
(except for the rightmost node, which corresponds to the last element of
$\sigma$), it follows 
from Proposition \ref{Rn_trees}
that the descents of a permutation $\sigma\in\R_n$
correspond to the leaves of $T(\sigma)$ minus the rightmost node. Thus for
$\sigma\in\R_n$ we have that $2^{d(\sigma)}$ is the number of leaves in
$T(\sigma)$ minus one, which equals 
the number of of nodes of $T(\sigma)$ with two
children. (The latter can be proved easily by induction.)

For a min-tree $T$ and a node $s$ of $T$ with two children, let $F_s(T)$ be
the tree obtained by switching the left and the right subtrees of $T$.
(This operation is called the \emph{Foata-Strehl action} on the permutation
encoded by $T$; see \cite{HetyeiReiner}.)
For example, if $T$ is the tree $T(423516)$ shown
above, then $F_2(T)$ is the tree shown in Figure \ref{sample_tree_2}.

\begin{figure}[!h]
\begin{center}
\input{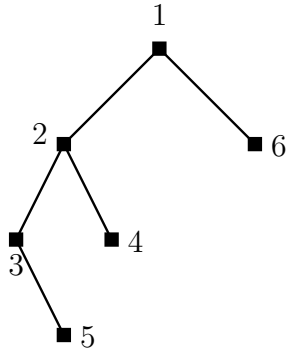}
\end{center}
\caption{The tree $F_2(T(423516))$}\label{sample_tree_2}
\end{figure}

Note that the action of $F_s$ preserves the set of nodes with two children and
does not create any nodes with a lone left child if the original tree
contained no such nodes. Hence the set $T(\R_n)$ 
is invariant under this action.
Observe also that the operators $F_s$ commute and satisfy $F_s^2 = 1$. Thus
these operators, viewed as operators on permutations corresponding to trees,
split the set $\R_n$ into orbits of size $2^{d(\sigma)}$, where~$\sigma$ is
any member of the orbit. It remains to show that each orbit contains exactly
one permutation in $\SSS_n$.

Given $\sigma\in\R_n$, there is a unique, up to order, 
sequence of operators $F_s$, where
$s$ is on the rightmost path, that, when applied to $T(\sigma)$, 
makes $n$ the rightmost node of the
resulting tree. An example is shown in Figure \ref{sample_tree_3}.
\begin{figure}[!h]
\begin{center}
\input{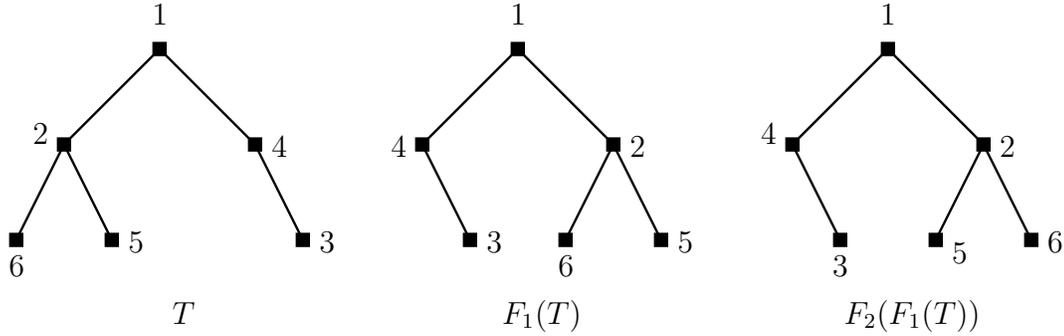}
\end{center}
\caption{The action of $F_1$ and $F_2$ on a min-tree}\label{sample_tree_3}
\end{figure}(One needs to find the closest ancestor of $n$ on the rightmost path and then
apply the corresponding operator to bring the node $n$ closer to the rightmost
path.) Once $n$ is the rightmost node, apply the operator $F_s$ to all nodes
$s$ with two children for which the condition of Proposition~\ref{simsun-trees}
is violated. We obtain a tree corresponding to a permutation
in~$\SSS_n$ in the orbit of $\sigma$. To see that each orbit contains only one
member of~$\SSS_n$, observe that the action of $F_s$ preserves the sequence of
elements on the path from $1$ to $k$ for each $k$, and given the sequence of
ancestors for each $k\in [n]$, there is a unique way of arranging the elements
of $[n]$ to form a min-tree satisfying the conditions of Proposition
\ref{simsun-trees}: first, set the path from $1$ to $n$ to be the rightmost
path, and then set all lone children to be right children, and for all nodes
with two children, set the greater element to be the left child.

The proof is now complete.
\end{proof}

Table \ref{alt_descent_set_cd} lists
the polynomials $\hat{\Phi}_n(\ccc,\ddd)$ for $n\leq 6$.
\begin{table}[h!]
$$
\begin{array}{l|l}
n\  &\ \hat{\Phi}_n(\ccc,\ddd)\\ \hline
1\  &\ 1 \\[5pt]
2\  &\ \ccc\\[5pt]
3\  &\ 2\ccc^2 - \ddd \\[5pt]
4\  &\ 5\ccc^3 - 2 (\ccc\ddd + \ddd\ccc) \\[5pt]
5\  &\ 16\ccc^4 - 7 (\ccc^2\ddd + \ddd\ccc^2)
                          - 5 \ccc\ddd\ccc + 4 \ddd^2 \\[5pt]
6\  &\ 61\ccc^5 - 26 (\ccc^3\ddd + \ddd\ccc^3)
                                   - 21 (\ccc\ddd\ccc^2 + \ccc^2\ddd\ccc)
                                    + 10 \ddd\ccc\ddd 
                                    + 12 (\ccc\ddd^2 + \ddd^2\ccc)
\end{array}
$$
\caption{The polynomials $\hat{\Phi}_n(\ccc,\ddd)$}
\label{alt_descent_set_cd}
\end{table}
\section{Shapiro-Woan-Getu permutations}
\label{SWG_permutations}

In this section we take a closer look at the class of permutations
which we denoted by $\mathcal{R}_n$ in Section \ref{cd-stuff}. 
Recall that $\mathcal{R}_n$ is the set of permutations with no consecutive
(double) descents and no descent at the end.
They appear
in the paper \cite{SWG} by Shapiro, Woan, and Getu, hence the section title,
who call them \emph{reduced} permutations. 
The paper studies enumeration of permutations by the number of runs
or slides, and in~\cite[Sec.~11.1]{PostnikovReinerWilliams}
Postnikov, Reiner, and Williams put these results in the context of
structural properties of permutohedra: for instance, the polynomial encoding
the distribution of permutations in $\mathcal{R}_n$ by the number of descents
is the $\gamma$-polynomial of the classical permutohedron.

In Section \ref{cd-stuff}, we found 
the number $R_n$ of SWG permutations of size $n$ to be the sum of 
absolute values of coefficients of a $(\ccc,\ddd)$-polynomial that, when
expanded in terms of $\aaa$ and $\bbb$, gave the generating function for the
alternating descent set statistic. Shapiro, Woan, and Getu provide a generating
function for $R_n$:
$$
R(x) := \sum_{n\geq 0} R_n\cdot {x^n\over n!} =
1 + {2\tan(x\sqrt{3}/2)\over \sqrt{3} - \tan(x\sqrt{3}/2)}
$$
(we put $R_0=1$). Observe that $R(x)R(-x) = 1$, a property that $R(x)$ shares
with $e^x$ and $\tan x + \sec x$, which are the two fundamental generating
functions in the analisys done in previous sections. There is a
further resemblance with the Euler numbers $E_n$ if one looks at the logarithm
of $R(x)$:
\begin{equation}\label{ln_Rx}
\ln(R(x)) =-x+2\sum_{n\geq 0}R_{2n}\cdot{x^{2n+1}\over{(2n+1)}!}.
\end{equation}
Comparing with
\begin{equation}\label{ln_tan_sec}
\ln(\tan x + \sec x) = \sum_{n\geq 0} E_{2n}\cdot {x^{2n+1}\over (2n+1)!},
\end{equation}
we see that taking the logarithm has a similar effect on both
$R(x)$ and $\tan x + \sec x$ of taking the even part and integrating,
except that for $R(x)$ all coefficients excluding that of $x$ are doubled.

The fact that
\begin{equation}\label{int_sec}
\int \sec x\ dx = \ln(\tan x + \sec x)
\end{equation}
(omitting the arbitrary constant of integration) has been used in the proof
of Theorem~\ref{f-eulerian}. This textbook integral formula can be proved
combinatorially using the exponential formula for generating functions
(see \cite[Sec.\ 5.1]{StanleyEC2}).
Given an up-down permutation $\sigma$,
divide $\sigma$ into blocks by the following procedure. Put the 
subword of $\sigma$ starting at the beginning of $\sigma$ and ending at 
the element equal to $1$ in the
first block, and remove this block from $\sigma$. In the resulting word,
find the \emph{maximum} element $m_2$ and put the subword consisting of initial
elements of the word up to, and including,~$m_2$ in the second block, and
remove the second block. In the remaining word, find the \emph{minimum} element
$m_3$, and repeat until there is nothing left, alternating between cutting
at the minimum and at the maximum element of the current word.
For example, for $\sigma = 5 9 3 4 1 8 6 7 2$, the blocks would be
$59341$, $8$, and $672$. Note that given the blocks one can uniquely
recover the order in which they must be concatenated to form the original
permutation $\sigma$. Indeed, the first block is the one containing $1$,
the second block contains the largest element not in the
the first block, the third block contains the smallest element not in the
first two blocks, and so on. Thus to construct an up-down permutation 
of size $n$ we need
to divide the elements of $[n]$ into blocks of \emph{odd} size, then determine
the order of concatenation using the above principle, and then arrange the
elements of odd numbered blocks in up-down order and those of even numbered
blocks in down-up order. There are $E_{k-1}$ ways to arrange the elements
in a block of size $k$ for odd $k$, and $0$ ways for even $k$ since we do not
allow blocks of even size. Thus (\ref{ln_tan_sec}),
which is equivalent to (\ref{int_sec}), follows from the exponential
formula. This argument ``combinatorializes'' the proof of
Theorem \ref{f-eulerian}.
It would be nice to give a similar argument for reduced permutations
$\mathcal{R}_n$.

\begin{problem}
Find a combinatorial proof of the formula (\ref{ln_Rx}) for $\ln(R(x))$.
\end{problem}

Another problem emerging from the results of
Section \ref{cd-stuff} is the following.

\begin{problem}
Give a combinatorial interpretation of the coefficients of the polynomial
$\hat{\Phi}_n(\ccc,\ -\ddd)$ by partitioning the set $\mathcal{R}_n$ into
classes corresponding to the $F_{n-1}$ monomials.
\end{problem}

It is worth pointing out here that even though one can split
$\mathcal{R}_n$ into $F_{n-1}$ classes corresponding to $(\ccc,\ddd)$-monomials
by descent set, like it was done for simsun permutations in Section
\ref{cd-stuff}, the resulting polynomial is different from
$\hat{\Phi}_n(\ccc,\ -\ddd)$. There are a few hints on what the correct
way to refine permutations in $\mathcal{R}_n$ could be. The coefficient
of $\ccc^{n-1}$ in $\hat{\Phi}_n(\ccc,\ -\ddd)$ is the Euler number $E_{n}$,
and the set $\mathcal{R}_n$ includes at least three kinds of permutations
mentioned in this paper that are counted by $E_{n}$: alternating
permutations ending with an ascent, simsun permutations, and
permutations $\sigma\in\S_{n}$ such that $\sigma \circ (n+1)$ has no
$3$-descents. Values of $\hat{\Phi}_n$ for small $n$, including those
listed in Table~\ref{alt_descent_set_cd},
present evidence that the common coefficient of $\ccc^{n-3}\ddd$
and $\ddd\ccc^{n-3}$ is $\hat{A}(n-1,2)$ (the number of permutations
of size $n-1$ with exactly one alternating descent).

\section{A $q$-analog of Euler numbers}\label{Euler-q-analog}

Let $\hat{A}_n(t,q)$ denote the bivariate polynomial of Theorem
\ref{main_bivariate_identity}:
$$
\hat{A}_n(t,q) := \sum_{\sigma\in\S_n} 
t^{\hat{d}(\sigma)} q^{\hat{\imath}(\sigma)}.
$$
Then the alternating Eulerian polynomial $\hat{A}_n(t)$ is just the
specialization $t\hat{A}_n(t,1)$. We also noted earlier (Corollary
\ref{hat_code_mahonian}) that
$$
\hat{A}_n(1,q) = [n]_q!,
$$
the classical $q$-analog of the factorial defined by $[n]_q!
:= [1]_q[2]_q\cdots [n]_q$, where $[i]_q := 1+q+q^2+\cdots+q^{i-1}$.
One can ask about other specializations of $\hat{A}_n(t,q)$, such as the ones
with $t$ or $q$ set to $0$. Clearly, we have $\hat{A}_n(t,0) = 1$ because
the only permutation~$\sigma\in\S_n$ for which $\hat{\imath}(\sigma) = 0$
also satisfies $\hat{d}(\sigma) = 0$. The case of $t=0$ is more curious and
is the subject of this section.

We have $\hat{d}(\sigma) = 0$ if and only if $\sigma$ is an up-down
permutation. Thus $\hat{A}_n(0,1) = E_n$, and the specialization
$\hat{A}_n(0,q)$ gives a $q$-analog of the Euler number $E_n$ with coefficients
encoding the distribution of the number of alternating inversions among
up-down permutations. The following lemma is key in understanding this
$q$-analog.

\begin{lemma}\label{alt_inversion_criterion}
For a permutation $\sigma\in\S_n$, let
$\hat{\code}(\sigma) = (\hat{c}_1,\hat{c}_2,\ldots,\hat{c}_{n-1})$.
Then $\sigma$ is up-down (resp.,\ down-up) if and only if
$\hat{c}_i + \hat{c}_{i+1} \leq n-1-i$ (resp.,\ 
$\hat{c}_i + \hat{c}_{i+1} \geq n-i$) for all $i$.
\end{lemma}

\begin{proof}
This fact is just a special case of Lemma \ref{hatD_from_hatcode}.
\end{proof}

For various reasons it is more convenient to study the distribution of
$\hat{\imath}$ on down-up, rather than up-down, permutations. The $q$-analog 
obtained
this way from down-up permutations is essentially equivalent to
$\hat{A}_n(0,q)$, the difference being the reverse order of coefficients and
a power of $q$ factor. 
It follows from Lemma \ref{alt_inversion_criterion} that for a down-up
permutation $\sigma\in\S_n$, we have
\begin{equation}\label{min_down_up_alt_inv}
\hat{\imath}(\sigma) \geq (n-1) + (n-3) + (n-5) + \cdots = 
\left\lfloor {n^2\over 4} \right\rfloor.
\end{equation}
Therefore let $\Alt_n$ be the set of down-up
permutations in $\S_n$, and define
$$
\hat{E}_n(q) := q^{-\lfloor n^2/4 \rfloor}
\sum_{\sigma\in\Alt_n} q^{\hat{\imath}(\sigma)}.
$$
The values of $\hat{E}_n(q)$ for small $n$ are given in 
Table \ref{hatEnq}.
\begin{table}[t!]
$$
\begin{array}{l|l}
n & \hat{E}_n(q)\\\hline
0, 1, 2 & 1\\[3pt]
3       & 1+q\\[3pt]
4       & 2+2q+q^2\\[3pt]
5       & 2+5q+5q^2+3q^3+q^4\\[3pt]
6       & 5+12q+16q^2+14q^3+9q^4+4q^5+q^6\\[3pt]
7       & 5+21q+42q^2+56q^3+56q^4+44q^5+28q^6+14q^7+5q^8+q^9
\end{array}
$$
\caption{The polynomials $\hat{E}_n(q)$ for $n\leq 7$}
\label{hatEnq}
\end{table}
We have the following facts about $\hat{E}_n(q)$.

\begin{proposition}\label{hatEq_facts}
(a) The polynomial $\hat{E}_n(q)$ is monic and has degree
$\left\lfloor {(n-1)^2\over 4}\right\rfloor$.

\vskip4pt
\noindent
(b) $\hat{A}_n(0,q) = q^{\lfloor (n-1)^2/4 \rfloor}\cdot\hat{E}_n(q^{-1})$.

\vskip4pt
\noindent
(c) $\hat{E}_n(0) = c_{\lfloor n/2 \rfloor}$, the $\lfloor n/2 \rfloor$-th
Catalan number.
\end{proposition}

\begin{proof}
(a) By Proposition \ref{hat_code}, the unique permutation
$\sigma\in\S_n$ with the maximum possible number of alternating inversions
is the one for which $\hat{\code}(\sigma) = (n-1,n-2,\ldots,1)$. By Lemma
\ref{alt_inversion_criterion}, or by simply realizing that
$\sigma = n \circ 1 \circ (n-1) \circ 2 \circ \cdots$, one can see that
$\sigma\in \Alt_n$. We have $\hat{\imath}(\sigma) = n(n-1)/2$, and thus
the degree of $\hat{E}_n(q)$ is $n(n-1)/2 - \lfloor n^2/4 \rfloor
= \lfloor (n-1)^2/4 \rfloor$.

\vskip4pt
(b) This identity is an algebraic restatement of an earlier observation.

\vskip4pt
(c) The constant term $\hat{E}_n(0)$ of $\hat{E}_n(q)$ is the number
of permutations $\sigma\in\Alt_n$ with exactly $\lfloor n^2/4 \rfloor$
alternating inversions. By (\ref{min_down_up_alt_inv}), these are precisely
the permutations in $\Alt_n$ satisfying $\hat{c}_{i} + \hat{c}_{i+1} = n-i$
for odd $i$. Let $\sigma\in\Alt_n$ be a permutation with this property.

For $j\geq 1$, we have $\hat{c}_{2j} \geq n-2j - \hat{c}_{2j+1}
= \hat{c}_{2j+2} - 1$. Thus $\hat{c}_2, \hat{c}_4, \ldots,
\hat{c}_{2\lfloor n/2 \rfloor}$ is a strictly decreasing sequence of
non-negative integers satisfying $\hat{c}_{2j} \leq n-2j$
(for convenience, let $\hat{c}_n = 0$).
Reversing the sequence and reducing the
$k$-th term by~$k-1$
for all $k$
yields a bijective correspondence with sequences of
$\lfloor n/2 \rfloor$ non-negative integers whose $k$-th term does not exceed
$k-1$, and it is well known that there are $c_{\lfloor n/2 \rfloor}$ such
sequences. Since $\hat{c}_{2j-1}$
is uniquely determined by $\hat{c}_{2j}$, it follows that there are
$c_{\lfloor n/2 \rfloor}$ permutations $\sigma\in\Alt_n$ with
$\lfloor n^2/4 \rfloor$
alternating inversions.
\end{proof}

It is curious to note that the permutations in $\Alt_n$ with
$\lfloor n^2/4 \rfloor$ alternating inversions can be characterized in terms
of pattern avoidance, so that Proposition \ref{hatEq_facts}(c) follows 
from a result of Mansour \cite{Mansour} stating that the number of 
$312$-avoiding down-up permutations of size $n$ is $c_{\lfloor n/2 \rfloor}$.

\begin{proposition}\label{alt_avoiding_312}
A permutation $\sigma\in\Alt_n$ has $\hat{\imath}(\sigma) = \lfloor
n^2/4 \rfloor$ if and only if $\sigma$ is $312$-avoiding.
\end{proposition}

The following lemma implies the above proposition and is useful in the later
discussion as well.

\begin{lemma}\label{alt_inversions_31-2}
For a permutation $\sigma = \sigma_1\sigma_2\cdots \sigma_n\in\Alt_n$,
the number $\hat{\imath}(\sigma)$ is
equal to $\lfloor n^2/4\rfloor$ plus the number of occurrences of
the generalized pattern $31$-$2$ (that is, the number of pairs of indices
$i < j$ such that $\sigma_{i+1} < \sigma_j < \sigma_{i}$).
\end{lemma}

\begin{proof}
For $i\in [n-1]$, define
$$
S_i := \left\{
\begin{array}{ll}
\{j\ |\ j>i\mbox{ and }\sigma_i>\sigma_j\}& \mbox{if $i$ is odd;}\\
\{j\ |\ j>i\mbox{ and }\sigma_i<\sigma_j\}& \mbox{if $i$ is even.}
\end{array}
\right.
$$
Thus $\hat{c}_i = |S_i|$. Let $i$ be \emph{odd}. Then 
$\sigma_i > \sigma_{i+1}$, so $i+1 \in S_i$ and for every $j>i+1$, either
$\sigma_j < \sigma_i$ or $\sigma_j > \sigma_{i+1}$, or both. Hence
$\{i+1,i+2,\ldots,n-1\}\subseteq S_i\cup S_{i+1}$ and $\hat{c}_i + \hat{c}_{i+1}
= n-1-i + |S_i\cap S_{i+1}|$. But $S_i\cap S_{i+1}$ is the set of indices
$j>i+1$ such that $\sigma_{i+1} < \sigma_j < \sigma_{i}$, i.e.\ the number
of occurrences of the pattern $31$-$2$ beginning at position $i$. Therefore
the total number of alternating inversions is $\sum_{i\ odd}\ (n-1-i) = 
\lfloor n^2/4 \rfloor$ plus the total number of occurrences of $31$-$2$. 
\end{proof}

\begin{proof_}{Proof of Proposition \ref{alt_avoiding_312}.}
Suppose that a permutation $\sigma\in\Alt_n$ has exactly
$\lfloor n^2/4 \rfloor$ alternating inversions
but is not $312$-avoiding. Choose a triple
$i < k < j$ such that $\sigma_k < \sigma_j < \sigma_i$ and the difference
$k-i$ is as small as possible. Suppose that $k-i \geq 2$. If
$\sigma_{k-1} < \sigma_j$, then we have $\sigma_{k-1} < \sigma_j < \sigma_i$,
contradicting the choice of $i$, $k$, and $j$. If $\sigma_{k-1} > \sigma_j$,
then we have $\sigma_k < \sigma_j < \sigma_{k-1}$, also contradicting the
choice of $i$, $k$, and $j$. Hence $k = i+1$, and we obtain a contradiction
by Lemma \ref{alt_inversions_31-2}.
\end{proof_}

In view of Lemma \ref{alt_inversions_31-2}, we can write $\hat{E}_n(q)$ as
$$
\hat{E}_n(q) = \sum_{\sigma\in\Alt_n} q^{\occp(\sigma)}
$$
where $\occp(\sigma)$ is the number of occurrences of $31$-$2$ in $\sigma$.
In what follows, we use this expression to show how a $q$-analog of a
combinatorial identity representing the Euler number $E_n$ as a weighted sum
of Dyck paths yields a refined identity of $\hat{E}_n(q)$.

First, we need to introduce Dyck paths, which are perhaps the most famous
combinatorial objects counted by Catalan numbers. A \emph{Dyck path} of
length $2m$ is a
continuous path consisting of line segments, or \emph{steps}, each of which
connects an integer point $(x,y)$ with either $(x+1,y-1)$ or $(x+1,y+1)$,
such that the path starts at~$(0,0)$, ends at $(2m,0)$, and never goes
below the $x$-axis, that is, contains no point with a negative $y$-coordinate.
The identity we are about to describe involves associating a certain weight
with every step of a Dyck path, defining the weight of the entire path to be
the product of the weights of the individual steps, and adding the weights of
all Dyck paths of length $2m$ to obtain $E_{2m}$ or $E_{2m+1}$, or, in the
case of the refined identity, $\hat{E}_{2m}(q)$ or $\hat{E}_{2m+1}(q)$.

For a step in a Dyck path, define the \emph{level} of that step to be the
$y$-coordinate of the highest point of the corresponding segment of the path.
Given a Dyck path $\mathcal{D}$ of length $2m$, let $\ell(i)$ be the level of
the $i$-th step of $\mathcal{D}$. Define
$$
w^e_{\mathcal{D},i}(q) := [\ell(i)]_q
$$
and
$$
w^o_{\mathcal{D},i}(q) := 
\left\{
\begin{array}{ll}
[\ell(i)]_q, & \mbox{if the $i$-th step is an up-step};\\

[\ell(i)+1]_q, & \mbox{if the $i$-th step is a down-step}.
\end{array}
\right.
$$
As mentioned above,
we set the weight of the entire path to be the product of step
weights:
\begin{eqnarray*}
w^e_{\mathcal{D}}(q) &=& \prod_{i=1}^{2m} w^e_{\mathcal{D},i}(q);\\
w^o_{\mathcal{D}}(q) &=& \prod_{i=1}^{2m} w^o_{\mathcal{D},i}(q). 
\end{eqnarray*}

\begin{theorem}\label{weighted_path_identity}
We have
$$
\sum_{\mathcal{D}} w^e_{\mathcal{D}}(q) =
\hat{E}_{2m}(q)
$$
and
$$
\sum_{\mathcal{D}} w^o_{\mathcal{D}}(q) =
\hat{E}_{2m+1}(q),
$$
where both sums are taken over all Dyck paths of length $2m$.
\end{theorem}

\begin{figure}[h!]
\begin{center}
\input{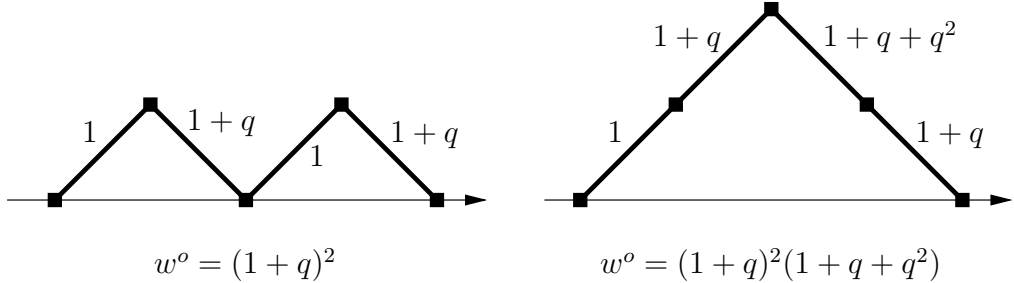}
\end{center}
\caption{Weighted Dyck paths adding up to $\hat{E}_5(q)$}
\label{paths2}
\end{figure}

For example, for $m=2$ there are two Dyck paths, shown in Figures
\ref{paths1} and \ref{paths2} with step weights given by
$w^e_{\mathcal{D},i}(q)$ and $w^o_{\mathcal{D},i}(q)$.
\begin{figure}[h!]
\begin{center}
\input{paths1.pstex_t}
\end{center}
\caption{Weighted Dyck paths adding up to $\hat{E}_4(q)$}
\label{paths1}
\end{figure}
%
From these weighted paths, we get
$$
1 + (1+q)^2 = 2 + 2q + q^2 = \hat{E}_4(q)
$$
and
$$
(1+q)^2 + (1+q)^2(1+q+q^2) = 2 + 5q + 5q^2 + 3q^3 + q^4 = \hat{E}_5(q).
$$

In the classical case $q=1$, the identities of Theorem
\ref{weighted_path_identity} are due to
Fran\c{c}on and Viennot \cite{FranconViennot},
and are discussed in a broader context in
the book 
\cite[Sec.\ 5.2]{GouldenJackson} by Goulden and Jackson.
The proof of our identities
is a refinement of the original argument.

\vskip12pt
\begin{proof_}{Proof of Theorem \ref{weighted_path_identity}.}
Fix a positive integer $n>1$, and let $m = \lfloor n/2 \rfloor$.
Recall that in Section \ref{cd-stuff} we associated to a permutation $\sigma
\in \S_n$ an increasing planar binary tree $T(\sigma)$ with vertex set $[n]$.
Extending the argument in the proof of Lemma \ref{no_consecutive_descents},
we conclude that $\sigma$ is in $\Alt_n$ if and only if the tree 
$T(\sigma)$ has no vertices with a lone child, except for the rightmost
vertex in the case of even $n$, which has a lone left child. For $\sigma
\in\Alt_n$, define the corresponding Dyck path $\mathcal{D}(\sigma)$ of
length $2m$ as follows: set the $i$-th step of the path to
be an up-step if vertex $i$ of $T(\sigma)$ has at least one child, and set the
$i$-th step to be a down-step if vertex $i$ is a leaf of $T(\sigma)$.\
We leave it as an exercise for the reader to check that $\mathcal{D}(\sigma)$
is a valid Dyck path.

Fix a Dyck path $\mathcal{D}$ of length
$2m$. We claim that
\begin{equation}\label{fixed_path_claim}
\sum_{\sigma\in\Alt_n\ :\ \mathcal{D}(\sigma)=\mathcal{D}} q^{\occp(\sigma)}
= \left\{\begin{array}{ll}
w^e_{\mathcal{D}}(q), & \mbox{if $n$ is even};\\
w^o_{\mathcal{D}}(q), & \mbox{if $n$ is odd}.
\end{array}\right.
\end{equation}
To prove the claim, consider for every $i$ the subtree $T_i(\sigma)$ obtained
from $T(\sigma)$ by removing all vertices labeled with numbers greater than
$i$. For the sake of clarity, one should imagine the ``incomplete'' tree
$T_i(\sigma)$ together with ``loose'' edges indicating those edges with
parent vertices in $T_i(\sigma)$ that appear when $T_i(\sigma)$ is completed
to~$T(\sigma)$. For even $n$ one should also think of a loose edge directed
to the right coming out
of the rightmost vertex of every tree $T_i(\sigma)$ including
$T_n(\sigma) = T(\sigma)$ --- this way the number of edges coming out of a
vertex of $T_i(\sigma)$ is always $0$ or $2$.

Observe that for $1\leq i \leq 2m$,
the number of loose edges of  $T_i(\sigma)$ is equal to $y_\mathcal{D}(i)+1$,
where $y_\mathcal{D}(i)$ is the $y$-coordinate of the
point of $\mathcal{D}$ whose $x$-coordinate is $i$.
Indeed, $T_1(\sigma)$ has
two loose edges, and $T_{i+1}(\sigma)$ is obtained from $T_i(\sigma)$ by
attaching a non-leaf to a loose edge, thus increasing the number of loose edges
by one, if the $i$-th step of $\mathcal{D}$ is an up-step, or by attaching
a leaf to a loose edge, thus reducing the number of loose edges by one, if
the $i$-th step is a down-step. Hence we can count the number of permutations
$\sigma\in\Alt_n$ with $\mathcal{D}(\sigma)=\mathcal{D}$ by multiplying
together the number of possibilities to attach a vertex labeled $i+1$ to
$T_i(\sigma)$ to form $T_{i+1}(\sigma)$ for all $1\leq i\leq n-1$.
The number of valid places to attach vertex $i+1$ is equal to the number
of loose edges in $T_i(\sigma)$ unless $i+1$ is a leaf of $T(\sigma)$ and
$n$ is even, in which case we have one fewer possibilities, because
we are not allowed to make the rightmost vertex a leaf.
Note that the level $\ell(i)$ of the $i$-th step of $\mathcal{D}$ is equal to
$y_\mathcal{D}(i)$ if it is an up-step, or $y_\mathcal{D}(i)+1$ if it is
a down-step.
Comparing with the choice of step weights, we conclude that
the number of possibilities to attach vertex $i+1$ is
$w^e_{\mathcal{D},{i+1}}(1)$
if $n$ is even, or $w^o_{\mathcal{D},{i+1}}(1)$ if $n$ is odd.
(For odd $n$ and $i=n-1$ the latter assertion makes no sense as $\mathcal{D}$
does not have an $n$-th step; however, there is just one way to attach the
last vertex, so the counting argument is not affected.)

The above computation proves the $q=1$ case of (\ref{fixed_path_claim}).
To prove the general claim, we need to show that if there are $p$ possibilities
to attach vertex $i+1$ to a loose edge of $T_i(\sigma)$, then the number of
occurrences of the $31$-$2$ pattern ``induced'' by the attachment is $0$
for one of the possibilities, $1$ for another possibility, $2$ for another,
and so on, up to $p-1$. Then choosing a place to attach vertex $i+1$ would
correspond to choosing a term from $1 + q + q^2 + \cdots +q^{p-1} =
[p]_q$, the weight of the $i$-th step of~$\mathcal{D}$, which is a factor in
the total weight of~$\mathcal{D}$, and (\ref{fixed_path_claim}) would follow.

It remains to specify which occurrences of $31$-$2$ in $\sigma$ are induced
by which vertex of~$T(\sigma)$. Suppose there are $p$ possible places to attach
vertex $i+1$. Order these places according to the topological order of tree
traversal, and suppose we choose to put vertex $i+1$ in the $k$-th place in
this order. Let $r_1$, $r_2$, \ldots, $r_{k-1}$ be the numbers of the vertices
immediately following the first $k-1$ places in the topological order, and
let~$a_j$ denote the label of the rightmost vertex of the eventual subtree of
$T(\sigma)$ rooted at what is currently the $j$-th of these $k-1$ places.
Although $a_j$ is not determined at the time vertex $i+1$ is attached, it is
certain that $r_j < i+1 <a_j$ and that $a_j$ and $r_j$ will be consecutive
elements of $\sigma$, with $i+1$ located somewhere to the right,
resulting in an occurrence of $31$-$2$. Thus the
choice to put vertex $i+1$ in the $k$-th available place induces $k-1$
occurrences of $31$-$2$, one for each $1\leq j \leq k-1$.
It is not hard to
check that each occurrence of $31$-$2$ is induced by some vertex of
$T(\sigma)$, namely, the vertex corresponding to the rightmost element forming
the pattern, in the way described above.

\begin{figure}[h!]
\begin{center}
\input{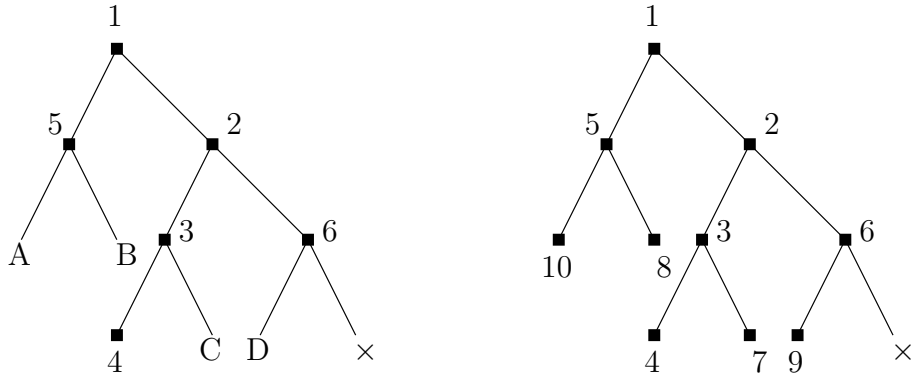}
\end{center}
\caption{An intermediate tree $T_6(\sigma)$ and
its completion $T(\sigma)$}
\label{tree_construction}
\end{figure}

Let us illustrate the argument with
an example. The left side of Figure \ref{tree_construction} shows
the tree $T_6(\sigma)$ for some $\sigma \in \Alt_{10}$, with the four potential
places for vertex $7$ marked A, B, C, and D. If vertex $7$ is put in position
A, then it induces no occurrences of $31$-$2$. If it is put in position B,
it induces one occurrence of $31$-$2$ as the triple $a5$-$7$ is created, where
$a$ stands for the number of the rightmost vertex in the subtree rooted at
A in the eventual tree. If vertex $7$ is put in position C, then in addition
to the triple $a5$-$7$, one obtains a second $31$-$2$ triple $b1$-$7$. Finally,
putting vertex $7$ in position D results in a third $31$-$2$ triple
$c2$-$7$. (Here $b$ and $c$ are defined by analogy with $a$.) On the right
side of Figure \ref{tree_construction} we have a possible completion of the
tree on the left, which corresponds to the permutation
$\sigma=10\ 5\ 8\ 1\ 4\ 3\ 7\ 2\ 9\ 6$.

The theorem now follows by taking the sum of (\ref{fixed_path_claim}) over
all Dyck paths $\mathcal{D}$ of length $2m$.
\end{proof_}

\end{document}